\documentclass{amsart}

\usepackage{latexsym,amsfonts,amssymb,exscale,enumerate,comment}
\usepackage{amsmath,amsthm,amscd}

\usepackage{url}
\usepackage[bookmarks=true,%
    colorlinks=true,%
    linkcolor=blue,%
    citecolor=blue,%
    filecolor=blue,%
    menucolor=blue,%
    urlcolor=blue,%
    breaklinks=true]{hyperref}

\input xy
\usepackage[all]{xy}
\xyoption{line}
\xyoption{arrow}
\xyoption{color}
\SelectTips{cm}{}

\usepackage{tikz}
\usetikzlibrary{shapes,snakes}
\usetikzlibrary{decorations.markings}
\usetikzlibrary{decorations.pathreplacing}
\tikzstyle directed=[postaction={decorate,decoration={markings,
    mark=at position #1 with {\arrow{>}}}}]

\newcommand{\hackcenter}[1]{
 \xy (0,0)*{#1}; \endxy}

\usepackage{tikz-cd}
\tikzset{->-/.style={decoration={
  markings,
  mark=at position #1 with {\arrow{>}}},postaction={decorate}}}

\tikzset{middlearrow/.style={
        decoration={markings,
            mark= at position 0.5 with {\arrow{#1}} ,
        },
        postaction={decorate}
    }
}

\usepackage{bm}

\def\u0{{\underline{0}}}

\def\la{\langle}
\def\ra{\rangle}

\newcommand{\onel}{\1_{\lambda}}

\def\P{\mathsf{P}}

\theoremstyle{plain}
\newtheorem{theorem}{Theorem}

\newtheorem{corollary}[theorem]{Corollary}
\newtheorem{proposition}[theorem]{Proposition}
\newtheorem{lemma}[theorem]{Lemma}

\newtheorem{notation}[theorem]{Notation}

\theoremstyle{definition}

\newtheorem{definition}[theorem]{Definition}

\theoremstyle{definition}
\newtheorem{remark}[theorem]{Remark}

\numberwithin{equation}{section}
\numberwithin{theorem}{section}

\newcommand{\maps}{\colon}
\newcommand{\und}[1]{\underline{#1}}
\newcommand{\oUcat}{\mf{U}}
\newcommand{\onen}{\1_{n}}
\newcommand{\posbubd}[2]{\mathord{
\begin{tikzpicture}[baseline = 0,scale=1.2]
  \draw[<-,thick,black] (0,0.4) to[out=180,in=90] (-.2,0.2);
  \draw[-,thick,black] (0.2,0.2) to[out=90,in=0] (0,.4);
 \draw[-,thick,black] (-.2,0.2) to[out=-90,in=180] (0,0);
  \draw[-,thick,black] (0,0) to[out=0,in=-90] (0.2,0.2);
 \node at (0,-.1) {$\scriptstyle{#1}$};
   \node at (0.3,0.5) {$\scriptstyle{\lambda}$};
   \node at (-0.2,0.2) {$\color{black}\bullet$};
   \node at (-0.4,0.2) {$\color{black}\scriptstyle{#2}$};
\end{tikzpicture}
}}
\newcommand{\negbubd}[2]{\mathord{
\begin{tikzpicture}[baseline = 0,scale=1.2]
  \draw[->,thick,black] (0.2,0.2) to[out=90,in=0] (0,.4);
  \draw[-,thick,black] (0,0.4) to[out=180,in=90] (-.2,0.2);
\draw[-,thick,black] (-.2,0.2) to[out=-90,in=180] (0,0);
  \draw[-,thick,black] (0,0) to[out=0,in=-90] (0.2,0.2);
 \node at (0,-.1) {$\scriptstyle{#1}$};
   \node at (-0.3,0.5) {$\scriptstyle{\lambda}$};
   \node at (0.2,0.2) {$\color{black}\bullet$};
   \node at (0.45,0.2) {$\color{black}\scriptstyle{#2}$};
\end{tikzpicture}
}}

\newcommand{\posbubdfffsl}[2]{\mathord{
\begin{tikzpicture}[baseline = 0,scale=1.2]
  \draw[<-,thick,black] (0,0.4) to[out=180,in=90] (-.2,0.2);
  \draw[-,thick,black] (0.2,0.2) to[out=90,in=0] (0,.4);
 \draw[-,thick,black] (-.2,0.2) to[out=-90,in=180] (0,0);
  \draw[-,thick,black] (0,0) to[out=0,in=-90] (0.2,0.2);
 \node at (0,-.1) {$\scriptstyle{#1}$};
   \node at (-0.2,0.2) {$\color{black}\bullet$};
   \node at (-0.55,0.2) {$\color{black}\scriptstyle{#2}$};
\end{tikzpicture}
}}
\newcommand{\negbubdfffsl}[2]{\mathord{
\begin{tikzpicture}[baseline = 0,scale=1.2]
  \draw[->,thick,black] (0.2,0.2) to[out=90,in=0] (0,.4);
  \draw[-,thick,black] (0,0.4) to[out=180,in=90] (-.2,0.2);
\draw[-,thick,black] (-.2,0.2) to[out=-90,in=180] (0,0);
  \draw[-,thick,black] (0,0) to[out=0,in=-90] (0.2,0.2);
 \node at (0,-.1) {$\scriptstyle{#1}$};
   \node at (0.2,0.2) {$\color{black}\bullet$};
   \node at (0.55,0.2) {$\color{black}\scriptstyle{#2}$};
\end{tikzpicture}
}}

\newcommand{\oddbubble}[1]{
\begin{tikzpicture}
\begin{scope} [ x = 10pt, y = 10pt, join = round, cap = round, thick, scale=3]
  \draw[-,thick,black] (0,0.3) to[out=180,in=90] (-.2,0.1);
  \draw[-,thick,black] (0.2,0.1) to[out=90,in=0] (0,.3);
 \draw[-,thick,black] (-.2,0.1) to[out=-90,in=180] (0,-0.1);
  \draw[-,thick,black] (0,-0.1) to[out=0,in=-90] (0.2,0.1);
  \draw[-,thick,black] (0.14,-0.05) to (-0.14,0.23);
  \draw[-,thick,black] (0.14,0.23) to (-0.14,-0.05);
 \node at (0,-.25) {$\scriptstyle{#1}$};
\end{scope}
\end{tikzpicture}}

\newcommand{\refequal}[1]{\xy {\ar@{=}^{#1}
(-1,0)*{};(1,0)*{}};
\endxy}

\renewcommand{\d}{{\sf d}}
\newcommand\NH{\mathsf{NH}}
\newcommand\ONH{\mathsf{ONH}}
\newcommand\OPol{\mathsf{OPol}}
\newcommand\OL{\mathsf{O}\Lambda}
\newcommand{\del}{\partial}
\newcommand{\xt}{\widetilde{x}}
\newcommand\symL[1]{\Lambda_{#1}}
\newcommand\Pol[1]{\mathsf{P}_{#1}}
\newcommand{\sch}{\mathfrak{s}}

\newcommand{\To}{\Rightarrow}

\newcommand{\Hom}{{\rm Hom}}

\renewcommand{\to}{\rightarrow}

\newcommand{\qrk}{\mathrm{rk}_q}

\def\ker{{\mathrm{ker}}}
\def\im{{\mathrm{im}}}
\def\grdim{{\mathrm{grdim}}}

\def\Id{\mathrm{Id}}

\def\mf{\mathfrak}

\def\Mod{{\mathrm{Mod}}}
\def\sMod{{\mathrm{sMod} }}

\def\sProj{{\mathrm{sProj} }}
\def\xMod{{\mathrm{Mod}_{X}}}
\def\xRep{{\mathrm{Rep}_{X}}}
\def\xProj{{\mathrm{Proj}_{X}}}

\newcommand{\qbin}[2]{
\left[
 \begin{array}{c}
 #1 \\
 #2 \\
 \end{array}
 \right]
}
\newcommand{\qbins}[2]{
\left[
 \begin{array}{c}
 \scs #1 \\
 \scs #2 \\
 \end{array}
 \right]
}

\numberwithin{equation}{section}
\def\LV#1{\textcolor[rgb]{0.00,0.49,0.25}{[LV: #1]}}%
\def\ME#1{\textcolor[rgb]{0.40,0.00,0.90}{[ME: #1]}}%

\let\tilde=\widetilde

\let\epsilon=\varepsilon

\usepackage{bbm}
\def\C{{\mathbb{C}}}
\def\N{{\mathbbm N}}

\def\Z{{\mathbbm Z}}

\def\cal#1{\mathcal{#1}}%
\def\1{\mathbbm{1}}%
\def\nn{\notag}

\def\flr#1{\left\lfloor{#1}\right\rfloor}
\newcommand{\funct}[4]{\left \{ \begin{array}{rcl} #1 & \rightarrow & #2 \\ #3 & \mapsto & #4 \end{array} \right.}

\def\la{\langle}
\def\ra{\rangle}

\renewcommand{\l}{\lambda}

\usepackage[normalem]{ulem}

\def\cal#1{\mathcal{#1}}

\newcommand\nc{\newcommand}
\nc\rnc{\renewcommand}
\nc\Kar{\operatorname{Kar}}
\nc\End{\operatorname{End}}

\newcommand{\scs}{\scriptstyle}

\nc\Sym{\operatorname{Sym}}

\allowdisplaybreaks

\title{Derived superequivalences for spin symmetric groups and odd $\mf{sl}_2$-categorifications}

\begin{document}
\setcounter{tocdepth}{3}

\author{Mark Ebert}
\email{markeber@usc.edu}
\address{Department of Mathematics\\ University of Southern California \\ Los Angeles, CA}

\author{Aaron D. Lauda}
\email{lauda@usc.edu}
\address{Department of Mathematics \& Department of Physics\\ University of Southern California \\ Los Angeles, CA}

\author{Laurent Vera}
\email{laurent.vera@rutgers.edu}
\address{Hill Center for the Mathematical Sciences \\ Rutgers University \\ Piscataway, NJ 08854-8019}
\date{\today}

\maketitle

\begin{abstract}
We show that actions of the odd categorification of $\mf{sl}_2$ induce derived superequivalences analogous to those introduced by Chuang and Rouquier.   Using Kang, Kashiwara, and Oh's action of the odd 2-category on blocks of the cyclotomic affine Hecke-Clifford algebra, our equivalences imply that blocks related by a certain affine Weyl group action are derived equivalent.  By recent results of Kleshchev and Livesey, we show this implies Brou\'{e}'s abelian defect conjecture for the modular representations of the spin symmetric group.
\end{abstract}

%
\section{Introduction}
%

Let $\Bbbk$ be an algebraically closed field of odd characteristic $p>0$.
Let $\Z_2 = \{ 1, z\}$ and consider the double cover $\tilde{\mf{S}}_n$ of the symmetric group $\mf{S}_n$
\[
\xymatrix{ 1 \ar[r] & \Z_2 \ar[r] & \tilde{\mf{S}}_n \ar[r] & \mf{S}_n \ar[r] & 1 }.
\]
%
%
Then there is a canonical central element $z\in \tilde{\mf{S}}_n$ and central idempotent $e_z:= (1-z)/2 \in \Bbbk \tilde{\mf{S}}_n$ giving an ideal decomposition
\[
\Bbbk \tilde{\mf{S}}_n = \Bbbk \tilde{\mf{S}}_n e_z \oplus \Bbbk \tilde{\mf{S}}_n(1-e_z),
\]
where the block $\Bbbk \tilde{\mf{S}}_n (1-e_z) \cong \Bbbk  \mf{S}_n$ correspond to the usual group algebra of the symmetric group.  The other block of  $\Bbbk \tilde{\mf{S}}_n$ given by $\Bbbk \tilde{\mf{S}}_n e_z$ are called \emph{spin blocks} of the symmetric group that we denote by $\Bbbk \mf{S}_n^-$. These blocks control the spin (or projective) representation theory of the symmetric group $\mf{S}_n$, see \cite{WW,Kl-book} and references therein  for an excellent introduction to the spin theory of $\mf{S}_n$.   Similarly, we denote by $\Bbbk \tilde{\mf{A}}_n$ the double covers of the alternating group $\mf{A}_n$.

The superblocks $\cal{B}^{(\rho,d)}$ of $\Bbbk \mf{S}_n^-$ are labelled by pairs $(\rho,d)$, where $\rho$ is a $\bar{p}$-core partition (see section~\cite[Section 2.3]{KLLi} and \cite{K1}) and $d$ is a non-negative integer such that $|\rho|+dp =n$.  It is known that the defect group $D^{(\rho,d)}$ of $\cal{B}^{(\rho,d)}$ is abelian when $d<p$.  Let $\mf{b}^{(\rho, d)}$ be the Brauer correspondent of the block $D^{(\rho,d)}$.  Assume then that $d<p$, then Brou\'{e}'s  defect conjecture for the spin symmetric group states that the block $\cal{B}^{(\rho,d)}$ with abelian defect group $D^{(\rho,d)}$ is derived equivalent to its Brauer correspondent $\mf{b}^{(\rho, d)}$.

As explained in \cite{KLLi}, unlike the blocks $\Bbbk \tilde{\mf{S}}_n (1-e_z) \cong \Bbbk  \mf{S}_n$ that are relatively well understood~\cite{CR,Marcus}, very  little is known about spin blocks of symmetric and alternating groups, and in particular, Brou\'{e}'s conjecture for them is wide open.  In this article, we prove the spin defect conjecture for symmetric and alternating groups by utilizing so-called \emph{odd} categorical $\mf{sl}_2$-actions.

\subsection{Chuang-Rouquier equivalences}

Since Chuang and Rouquier's pioneering work~\cite{CR}, the construction of highly nontrivial derived equivalences has been one of the most powerful applications arising from the theory of categorified quantum groups.  Their work showed that categorical $\mf{sl}_2$-actions gave rise to derived equivalences coming from a lift of the Weyl group action.  Here we will call these CR-equivalences, or Chuang-Rouquier derived equivalences.    They used these CR-equivalence to prove Brou\'{e}'s abelian defect conjecture for the modular representation theory of the symmetric group.  Later, Marcus used these results to prove Brou\'{e}'s conjecture for alternating groups~\cite{Marcus}.

Since Chuang and Rouquier's work, a full theory of categorified quantum groups has emerged~\cite{Lau1,KL1,KL3,Rou2,RouQH} including the discovery of 2-categories $\cal{U}(\mf{g})$ categorifying integral forms of symmetrizable quantum Kac-Moody algebras.
CR-equivalences were then
 extended by Cautis, Kamnitzer, and Licata~\cite{CautisKam} to more general Lie type $\mf{g}$, where they correspond to simple transpositions giving rise to categorical braid group actions of type $\mf{g}$.
%
Such actions can be understood as giving the braidings in categorifications of the Jones and HOMFLYPT polynomials by Khovanov and Khovanov-Rozansky homology~\cite{CautisKam,LQR,QR,MacWeb} and even the braiding in bordered approaches to Knot Floer homology ~\cite{ManionTrivalent}.  These equivalences have also been used to construct derived equivalences related to stratified flops in birational geometry~\cite{CKL,CautisKam-coherent,Cautis}.

\subsection{Odd CR-equivalences}
In this article, we construct new derived superequivalences coming from the `odd' categorification of $\mf{sl}_2$. This odd theory arose
%
as an attempt to provide a higher representation theoretic explanation for a phenomena discovered in link homology theory where Khovanov homology exhibited two distinct even and odd variants~\cite{ORS,Putyra}.
Ellis, Khovanov, and the second author initiated a program~\cite{EKL} to define odd analogs of quantum $\mf{sl}_2$ and related structures.  The result was the discovery of odd, noncommutative, analogs of many of the structures appearing in connection with $\mf{sl}_2$ categorification including odd analogs of the Hopf algebra of symmetric functions~\cite{EK,EKL}, the nilHecke algebra, and cohomologies of Grassmannians~\cite{EKL} and Springer varieties~\cite{LauR}.  Subsequent work has shown these odd categorifications extend to arc algebras and constructions of odd Khovanov homology for tangles~\cite{NV,NK20,OddArc}.

These investigations into odd categorification turned out to be closely connected with independent parallel investigations into  Kac-Moody superalgebra categorifications~\cite{KKT, KKO, KKO2}, with the odd categorification of $\mf{sl}_2$ lifting the rank one  Kac-Moody superalgebra.  These odd categorifications 
give categorifications of
the theory of covering Kac-Moody algebras~\cite{HillWang,ClarkWang, CHW,CHW2}.  Covering algebras $U_{q,\pi}(\mf{g})$ generalize quantum enveloping algebras, depending on an additional parameter $\pi$ with $\pi^2 = 1$.  When $\pi=1$, it reduces to the usual quantum enveloping algebra $U_{q}(\mf{g})$, while the $\pi = -1$ specialization recovers the quantum group of a super Kac-Moody algebra.

In the rank one case, the $\pi=1$ specialization is $U_{q}(\mf{sl}_2)$, while for $\pi=-1$ it gives the quantum group $U_{q}(\mf{osp}(1|2))$ associated with the superalgebra $\mf{osp}(1|2)$.  Following a categorification of the positive parts of these algebra in \cite{HillWang}, Ellis and the second author categorified the full rank one covering algebra proving a conjecture from \cite{ClarkWang}.  In doing so, a 2-supercategory $\mf{U}:=\mf{U}(\mf{sl}_2)$ was defined~\cite{Lau-odd} for the rank one covering algebra whose Grothendieck group recovers $U_{q,\pi}(\mf{sl}_2)$.  This categorification was later greatly simplified by Ellis and Brundan~\cite{BE2}, where the 2-supercategory formalism was better developed, building off of their work~\cite{BE1}.

Despite the discovery of an `odd categorification of $\mf{sl}_2$' over a decade ago, the analog of CR-equivalences has remained elusive.  While it is relatively straightforward to define analogous complexes to the CR-equivalences, proving invertibility has proven more challenging.   By utilizing new strategies for proving invertibility developed by the third author, we are able achieve the culmination of the odd categorification program and define odd analogs of the CR-complexes $\Theta$.   In doing so, we achieve stronger  results than those developed by Chuang and Rouquier in \cite{CR}.  We prove results analogous to those obtained by Rouquier~\cite{Rou2}, showing that the odd CR-complexes are invertible  in the additive framework of any integrable 2-representation. %

\subsection{Applications to spin defect conjecture and affine Hecke-Clifford algebras}
Brundan and Kleshchev discovered a link between blocks of the affine Clifford-Hecke superalgebra (or affine Sergeev algebra) and quantum Kac-Moody superalgebras~\cite{BK-Hecke}, see also \cite{Tsu3,Klesh-book}.
In much the same way that KLR-algebras (quiver Hecke algebras) in type $A$ can be viewed as a graded analog of the affine Hecke algebra via the Brundan-Kleshchev-Rouquier isomorphism~~\cite{BK1,Rou2}, Kang, Kashiwara, and Tsuchioka introduced quiver Hecke superalgebras that serve as graded analogs of the affine Hecke-Clifford and affine Sergeev superalgebras.

Let $I=I_{{\rm even}} \sqcup I_{{\rm odd}}$ and $A=(a_{ij})_{i,j \in I}$ a symmetrizable Cartan matrix satisfying various compatibility conditions to define a symmetrizable Kac-Moody superalgebra $\mf{g}$, see Section~\ref{sec:QHA}.  To this Cartan data and an additional choice of certain skew polynomials $\mf{Q}$,
Kang, Kashiwara, Tsuchioka introduced ~\cite{KKT} quiver Hecke superalgebras $R_n=R_n(\mf{Q})$ that categorify positive halves of quantum Kac-Moody superalgebras~$\mf{g}$ associated to the data~\cite{KKO,KKO2}.   To a dominant integral weight $\Lambda \in P_+$, these superalgebras admit cyclotomic quotients $R_n^{\Lambda}$ giving categorifications of the highest weight representations $V(\Lambda)$ of the quantum Kac-Moody superalgebra $\mf{g}$.

The quiver Hecke superalgebras can be extended to 2-supercategories $\mf{U}(\mf{g})$ introduced by Brundan and Ellis~\cite{BE2} that act via 2-representations on categories of modules over cyclotomic quiver Hecke superalgebras $R_n^{\Lambda}$.  For each $i\in I_{{\rm odd}}$, the 2-category contains a copy of the odd $\mf{sl}_2$ 2-superategory $\mf{U}(\mf{sl}_2)$ and for $i\in I_{{\rm even}}$ it contains the usual even 2-category $\cal{U}(\mf{sl}_2)$ from \cite{Lau1}.  Hence, in any 2-representation of $\mf{U}(\mf{g})$ and any $i \in I$, there is a derived superequivalence $\Theta_i$ coming from either our newly defined odd CR-complexes if $i$ is odd, or the original CR-complexes if $i$ is even.  In particular,   our results imply derived superequivalences between blocks $R_{\alpha}^{\Lambda}$ of cyclotomic quiver Hecke algebras related by an action of corresponding Weyl group $W=W_{\mf{g}}$ (see Corollary~\ref{cor:QHSAblocks}).

There is a weak Morita superquivalence relating   $R_n$ and its cyclotomic quotients $R_n^{\Lambda}$  to quiver Hecke-Clifford superalgebras $RC_n$ and its cyclotomic quotients $RC_n^{\Lambda}$, also introduced in \cite{KKT}.   Kang, Kashiwara, and Tsuchioka showed that quiver Hecke-Clifford superalgebras admit  an
  analog of the Brundan-Kleshchev-Rouquier isomorphism, giving isomorphisms between blocks of cyclotomic affine Hecke-Clifford superalgebras $\cal{ACH}_{\beta}^{\Lambda}$ (or cyclotomic affine Sergeev superalgberas $\overline{\cal{ACH}}_{\beta}^{\Lambda}$ in the degenerate case) with blocks of their newly introduced quiver Hecke-Clifford algebras $RC_{\beta}^{\Lambda}$.  Hence, our derived equivalences for blocks of $R_{\alpha}^{\Lambda}$ give rise to related equivalences for blocks of affine Hecke-Clifford and affine Sergeev superalgebras.

The spin symmetric group $\Bbbk\mf{S}_n^-$ superalgebra of order $n$ appears in connection with the level one cyclotomic Sergeev superalgebra $\overline{ACH}_n^{\Lambda_0}$.  In particular, there is an isomorphism~\cite{Sergeev,Yam}
\begin{equation}
  \overline{ACH}_n^{\Lambda_0} \cong \Bbbk\mf{S}_n^- \otimes \mf{C}_n
\end{equation}
where $\mf{C}_n$ is the Clifford superalgebra.
Building off of foundational work of Kleshchev and Livesey~\cite{KLLi}, and categorical actions introduced in \cite{KKO,KKO2}, our derived equivalences for blocks of cyclotomic quiver Hecke superalgebras for $\Lambda=\Lambda_0$ implies the abelian defect conjecture for spin representations of the symmetric and alternating group (Theorem~\ref{thm:spindef}).

 \subsection{Specifics}
After recalling some basic notions in super theory in Section~\ref{sec:super},   Section \ref{sec:oddnil} introduces deformed cyclotomic quotients of the odd nilHecke algebras and shows they are Morita equivalent to a ring that can be thought of as an odd analog of the equivariant cohomology ring of a Grassmannian (see Theorem~\ref{thm:matrixMnk}).  Indeed, when coefficients are reduced modulo 2, this noncommuative ring reduces to the usual $GL(N)$-equivariant cohomology ring of a Grassmannian.

While writing this article, we became aware that Brundan and Kleshchev were independently working on the same problem with significant overlap to the results obtained here.  In particular, they have
independently constructed the Morita equivalence between deformed cyclotomic quotients and the odd equivariant cohomologies of Grassmannians~\cite{BK-new}.  They also construct the odd CR-complexes $\Theta$ and give an independent proof of the invertibility of the complex $\Theta$ acting on abelian 2-representations~\cite{BK-new}.  Our proofs of invertibility are entirely complementary.  They have defined 2-representations of $\mf{U}$ directly on these deformed Grassmannian superalgebras, analogous to the equivariant 2-representations from \cite{Lau2} defined in the even case.  They have significantly expanded the odd theory of symmetric functions and introduced bimodules that can be interpreted as odd equivariant cohomologies of two step flag varieties.    Their proof of invertibility of $\Theta$ then follows Chuang and Rouquier developing odd analogs of the theory of locally finite abelian 2-representations.

Here we work with the deformed cyclotomic quotients directly and extend to work of Kang, Kashiwara, and Oh to show that these deformed cyclotomic quotients admit 2-representations of $\mf{U}$, see Section~\ref{sec:2reps-def}.  This extension enables us to avoid quotienting $\mf{U}$ by the `odd bubble' (see \eqref{eq:def-odd-bubble}), whereas in \cite{BK-new} the odd bubble is set to zero.  While our setting is more general, this has no material effect in constructing derived equivalences for quiver Hecke superalgebras. (These two representations appear to factor through the quotient by the odd bubble.)

By relating the 2-representation on deformed cyclotomic quotients to a universal 2-representations constructed as a quotient of the 2-category $\mf{U}$ itself, we are able to define odd analogs   of the \emph{simple 2-representations} analogous to those studied by Rouquier in the even case~\cite{Rou2}.  These 2-representations are often more convenient than the minimal categorifications developed in \cite{CR} as they are additive rather than abelian.  Building off of the techniques developed by the third author, we can then study natural odd analogs of the Chuang-Rouquier complexes in Section~\ref{sec:def-superequiv} and prove in Theorem~\ref{thm:main} that these are invertible complexes in any integrable 2-representation.

\subsection{Acknowledgements}
 The authors are grateful to Raphael Rouquier and Matt Hogancamp for helpful conversations.  They are also grateful to Jon Brundan and Sasha Kleshchev for discussing their work in progress and clarifying aspects of Kleshchev-Livesey~\cite{KLLi}. The second author first became aware of the connection between the spin symmetric group and odd categorifications through conversations with Jon and Sasha. We are also grateful to Jon for pointing out a gap in an earlier version of this article.
A.D.L. and M.E were partially supported by NSF grants DMS-1902092 and DMS-2200419, the Army Research Office award W911NF-20-1-0075, and the Simons Foundation collaboration grant on New Structures in Low-dimensional Topology.

%
\section{Super theory} \label{sec:super}
%

\subsection{Superspaces}
Let $\Bbbk$ be a field with characteristic not equal to 2.   A \emph{superspace} is a $\Z_2$-graded vector space
$
V = V_{ \bar{0}} \oplus V_{ \bar{1}}.
$
For a homogeneous element $v \in V$, write $|v|$ for the parity of~$v$.  An \emph{even} linear map is a parity preserving $\Bbbk$-module map.  The usual tensor product of $\Bbbk$-vector spaces is again a superspace with
$
(V\otimes W)_{ \bar{0}} = V_{ \bar{0}} \otimes W_{ \bar{0}}\oplus V_{ \bar{1}} \otimes W_{ \bar{1}}
$
and
$
(V\otimes W)_{ \bar{1}} = V_{ \bar{0}} \otimes W_{ \bar{1}}\oplus V_{ \bar{1}} \otimes W_{ \bar{0}}.
$
Likewise, the tensor product $f\otimes g$ of two linear maps between superspaces is defined by
\begin{equation}
  (f\otimes g) (v \otimes w) := (-1)^{|g||v|} f(v) \otimes g(w).
\end{equation}
We write ${\rm SVect}$ for the category of superspaces and all linear maps. The underlying category $\und{{\rm SVect}}$ consisting of super spaces and even linear maps is a symmetric monoidal category with symmetric braiding $u \otimes v \mapsto (-1)^{|u||v|}v \otimes u$.
On a superspace, we have an involution
\begin{equation} \label{eq:varphi-inv}
  \varphi(y) = (-1)^{\vert y \vert}y .
\end{equation}

In this article, we will be primarily interested in supervector spaces equipped with an additional $\Z$-grading
\begin{equation} \label{eq:graded-super}
 V = \bigoplus_{n \in \Z}\bigoplus_{a \in \Z_2} V_{n,a}
\end{equation}

\subsubsection{$(q,\pi)$-integers}
Let $\cal{A}_{q,\pi}$ denote the ring $\Z[q,q^{-1},\pi]/(\pi^2-1)$.  In what follows, it is helpful to define $(q,\pi)$ quantum integers
 by setting
\begin{equation}
[n]_{q,\pi} := \frac{q^n - (\pi q)^{-n}}{q-(\pi q)^{-1}}
=
\left\{
  \begin{array}{ll}
    q^{n-1}+\pi q^{n-3} + \dots + \pi^{n-1}q^{1-n}, & \hbox{if $n\geq 0$, } \\
    -\pi^n\left(q^{-n-1} + \pi q^{-n-3} + \dots + \pi^{-n-1}q^{1+n} \right), & \hbox{if $n\leq 0$.}
  \end{array}
\right.
\end{equation}
for $n\in \Z$.
We then have analogs of the quantum factorials and quantum binomial coefficients
\begin{equation}
 [n]_{q,\pi}^! := [n]_{q,\pi} [n-1]_{q,\pi} \dots [1]_{q,\pi}, \quad
 \qbin{n}{r}_{q,\pi} := \frac{[n]^!_{q,\pi}}{[r]^!_{q,\pi} [n-r]^!_{q,\pi}}
\end{equation}
 Given a graded supervector space $V$ as in \eqref{eq:graded-super}, we  write
\[
\dim_{q,\pi} V :=  \sum_{n\in\Z} \sum_{a \in \Z_2} q^n\pi^a\dim V_{n,a}
\]
In what follows, we often just need the graded dimension $\grdim:= \dim_{q,\pi=1}$.    When working over $\Z$, rather than a field we write ${\rm rk}_{q,\pi}$ or ${\rm rk}_q$ for the corresponding ranks of free graded super or graded $\Z$-modules.

\subsection{Supercategories}

There are several variants of the definition of super monoidal category and 2-supercategory.  Here we follow the notation of \cite{BE1}, where a comparison to the notation from \cite{KKT,KKO} is provided at the end of the introduction in \cite{BE1}.

Supercategories, superfunctors, and supernatural transformations are all defined via enriched category theory over the symmetric monoidal category $\und{{\rm SVect}}$, see \cite[Definition 1.1]{BE1} for an unpacking of these definitions.  In particular, the Hom spaces $\Hom_{\cal{A}}(X,Y)$ of a supercategory $\cal{A}$ are equipped with superspace structures and composition and identities are given by even linear maps; a superfunctor $F \maps \cal{A} \to \cal{B}$ gives an even map of superspaces $\Hom_{\cal{A}}(X,Y) \to \Hom_{\cal{B}}(FX,FY)$ for all objects $X$, $Y$ of $\cal{A}$.

The notion of \emph{2-supercategory} can be defined in a similar way, by enriching over the category ${\rm SCat}$ of small supercategories and superfunctors with monoidal structure as in \cite[Definition 1.2]{BE2}.  In particular, a 2-supercategory $\mf{U}$ has for each pair of objects $\l,\mu$ a supercategory of morphisms $\Hom_{\mf{U}}(\lambda, \mu)$ and compositions and identities are given by monoidal superfunctors.  For our purposes, the key feature of 2-supercategories is that the interchange law must be replaced by a \emph{superinterchange law} so that
given 2-morphisms $u \maps p \To q \maps \lambda  \to \mu$, $u' \maps p'\To q' \maps \mu \to \nu$,  $v \maps q \To r \maps \lambda  \to \mu$, $v' \maps q'\To r' \maps \mu \to \nu$,  we have
\begin{equation} \label{eq:superinterchange}
    (u   u')\circ(v  v')
    \;\; = \;\; (-1)^{|u'| |v|}
    (u \circ v)  (u'\circ v'),
\end{equation}
where we denote horizontal composition by juxtaposition and vertical composition with $\circ$.

A graded $(Q, \Pi)$-supercategory is a graded supercategory $\cal{A}$ together with superfunctors
\[Q,Q^{-1},\Pi\maps \cal{A} \to \cal{A},\]
 an odd supernatural isomorphism $\zeta \maps \Pi \To I$ that is homogeneous of degree 0, and even supernatural isomorphisms $\sigma \maps Q \To I$, $\bar{\sigma} \maps Q^{-1} \To I$ that are homogeneous of degrees 1 and -1, respectively.  This data makes $Q$ and $Q^{-1}$ mutually inverse graded superequivalences and $\Pi$ a self-inverse graded superequivalence.  A graded $(Q,\Pi)$-2-supercategory can be defined similarly, where each hom category has the structure of a graded $(Q,\Pi)$-supercategory, see \cite[Definition 6.5]{BE1} for a precise definition.

\begin{notation}  \label{notation:direct}
 Given  $f=\sum_{j\in \Z}\sum_{a \in \Z_2}f_{j,a}q^j\pi^a \in \cal{A}_{q,\pi}$  and object $M$ of a  graded $(Q, \Pi)$-category, we write $fM$ to denote the direct sum
$\bigoplus_j(Q^{-j}\Pi^a M)^{\oplus f_{j,a}}$.
\end{notation}

\subsection{Superalgebras and supermodules} \label{sec:superalgebras}

A \textit{superalgebra} $A$ over the field $\Bbbk$ is a superspace together with a structure of $\Bbbk$-algebra such that $xy \in A_{i+j}$ for all $x\in A_{i}$ and $y\in A_{j}$. A superalgebra $A$ is said to be \textit{supercommutative} if $xy=(-1)^{ij}yx$ for all $x\in A_{i}$ and $y \in A_{j}$. Given superalgebras $A$ and $B$, the multiplication on the tensor product $A \otimes B$ is given on $\Z_2$ homogeneous elements by
\begin{equation} \label{eq:superalg-tensor}
 (a_1 \otimes b_1)(a_2 \otimes b_2) = (-1)^{|a_2||b_1|} (a_1a_2) \otimes (b_1b_2).
\end{equation}

A left (resp. right) supermodule over the superalgebra $A$ is a superspace $M$ together with a structure of left (resp. right) $A$-module such that $rm \in M_{i+j}$ (resp. $mr \in M_{i+j}$) for all $r\in A_{i}$ and $m \in M_{j}$.  We denote the supercategory of left $A$-supermodules by $\sMod(A)$.
The category $\Mod(A)$ of ordinary $A$-modules can be endowed with a $\Pi$-category structure, see \cite[Section 2.2]{KKT}.

The supercategory $\sMod(A)$ can be further equipped with the structure of a $\Pi$-supercategory.
If $M$ is a left (resp. right) supermodule, we define a left (resp. right) supermodule
	\[
		\Pi M = \left \{ \pi(m), \, m \in M \right\}
	\]
by $(\Pi M)_{i} = \left \{ \pi (m), \, m \in M_{i+\bar1} \right\}$ and $a\pi(m) = \pi (\varphi(a)m)$ (resp. $\pi(m)a = \pi(ma)$).
An even morphism of left (resp. right) supermodules is an even linear map $f$ such that $f(ay)=af(y)$. An odd morphism is an odd linear map $f$ such that $f(ay) = \varphi(a)f(y)$ (resp. $f(ya)=f(y)a$). Note that an odd morphism $M \rightarrow N$ is the same thing as an even morphism $M \rightarrow \Pi N$.  Then the identity map on the underlying vector space defines an odd supermodule isomorphism $\zeta_M\maps \Pi M \to M$ giving the category $\sMod(A)$ of $A$-supermodules and even $A$-linear homomorphisms the structure of a supercategory in the sense of \cite[Definition 2.1]{KKT}, or a $\Pi$-supercategory in the sense of \cite[Definition 1.7]{BE1} (see for example \cite[Example 1.8]{BE1}).

If in addition, $A$ is a $\Z$-graded superalgebra $A=\bigoplus_{i\in \Z}A_{i,\bar{0}} \oplus \bigoplus_{i\in \Z}A_{i,\bar{1}}$, we can consider the super category of graded supermodules.  This supercategory acquires a $(Q,\Pi)$-supercategory structure with $\Pi$ as above and $Q$ the grading shift functor defined on a $\Z$-graded supermodule $M=\bigoplus_{i\in \Z, j \in \Z_2}M_{i,j}$ by $(QM)_{i,j}=M_{i-1,j}$.

Define the super opposite $A^{{\rm sop}}$ of the superalgebra $A$ as $A^{{\rm sop}}_i := \{ a^{{\rm sop}} \mid a \in A_{i}, \}$ for $i \in \Z_2$ with multiplication
\[
a^{{\rm sop}} b^{{\rm sop}} = (-1)^{|a||b|}(ba)^{{\rm sop}}.
\]
Any left $A$-module can be regarded as a right $A^{{\rm sop}}$-module with $xa^{{\rm sop}} := (-1)^{|a||x|}a x$ for $\Z_2$-homogenous elements $a\in A$ and $x \in M$.

Given superalgebras $A$ and $B$, a superbimodule $M$ is a superspace $M$ and an even linear map $m_V\maps A \otimes V \otimes B \to M$ making $M$ into an $(A,B)$-bimodule in the usual sense.  A superbimodule homomorphism $f\maps M \to N$ is a linear map such that $f(avb)=(-1)^{|f||a|}af(v)b$. When $A$ and $B$ are graded superalgebras, the supercategory of graded $(A,B)$-superbimodules is also a $(Q,\Pi)$-supercategory~\cite[Example 1.8]{BE1}.  Given a graded $(A,B)$-superbimodule $M$, we have a functor $F_M \maps \sMod(B) \to \sMod(A)$ given by $N \mapsto M \otimes_{B} N$.  This functor can be equipped with the structure of a $(Q,\Pi)$-superfunctor~\cite[Section 7.5]{KKO2}.

 In what follows, we sometimes simplify terminology using algebra, module, and bimodule in place of superalgebra, supermodule, and super bimodule when it is clear from the context.

%
\section{Odd nilHecke algebra and deformed cyclotomic quotients}\label{sec:oddnil}
%

\subsection{The odd nilHecke algebra}
We recall here the odd nilHecke algebra introduced in \cite{EKL,KKT}.  This algebra is closely  related to the spin Hecke superalgebra appeared in earlier work of Wang~\cite{Wang};   many of the essential features of the odd nilHecke algebra including skew-polynomials appears much earlier in this and related works on spin symmetric groups~\cite{KW1,KW2,KW4}.
\begin{definition} \label{def-oddnil}
The \text{odd nilHecke algebra} $\ONH_n$ is the $\Z$-graded unital associative superalgebra generated by elements $x_1,\ldots,x_n$ of degree 2 and parity $\bar{1}$ and elements $\tau_1,\ldots,\tau_{n-1}$ of degree $-2$ and parity $\bar{1}$, subject to the relations
\begin{eqnarray}
& & \tau_i^2 = 0 , \quad \tau_i \tau_{i+1}\tau_i =
 \tau_{i+1}\tau_i \tau_{i+1},\\
 & &  x_i \tau_i + \tau_i  x_{i+1} =1, \quad
   \tau_i x_i + x_{i+1}\tau_i = 1,\\
& & x_i x_j + x_j x_i =0 \quad (i\neq j),  \quad
\tau_i \tau_j + \tau_j \tau_i =0 \quad (|i-j|>1),  \\
& & x_i \tau_j +\tau_j x_i = 0 \quad (i\neq j,j+1).
\end{eqnarray}
\end{definition}
For $w\in \mf{S}_n$ and a choice  of a reduced expression $w=s_{i_1}\cdots s_{i_\ell}$ in terms of simple transpositions $s_i=(i\quad i+1)$, define
$
\tau_w=\tau_{i_1}\cdots\tau_{i_\ell}.
$
Note that $\tau_w$ only depends on the reduced expression up to an overall sign.  For $w_0$ the longest word in $\mf{S}_n$ we fix a preferred choice of reduced expression
\begin{equation} \label{eq:reduced}
\tau_{w_0} = \tau_1(\tau_2\tau_1) \dots (\tau_{n-1} \dots \tau_1).
\end{equation}
One can show that
$
\tau_{w_0} = \tau_{n-1} (\tau_{n-2}\tau_{n-1}) \dots (\tau_1\dots\tau_{n-1} )
$, see \cite[(3.51)]{EKL}.
The elements
\begin{equation} \label{eq:idempotent}
  e_n :=   (-1)^{\binom{n}{3}}x_1^{n-1} x_2^{n-1} \dots x_n^0\tau_{w_0}, \\
\end{equation}
\begin{equation} \label{eq:idempotent'}
	e_n' :=   (-1)^{\binom{n}{3}}\tau_{w_0}x_n^{n-1} x_{n-1}^{n-2} \dots x_1^0,\\
\end{equation}
are idempotents of $\ONH_n$ and the left module $\ONH_ne_n$ is, up to grading shifts, the unique indecomposable projective $\ONH_n$-module, see \cite{EKL}.

 Fix a reduced expression $w=i_1\cdots i_r$ for each $w\in \mf{S}_n$.   If $\alpha=(\alpha_1,\ldots,\alpha_n)$ is an $n$-tuple, write $x^\alpha$ for $x_1^{\alpha_1}\cdots x_n^{\alpha_n}$. The basic properties of $\ONH_n$ are as follows.
\begin{proposition}[\cite{EKL}] \label{prop:oNH-facts} \hfill
\begin{enumerate}
\item Changing the choice of reduced expression only changes $\tau_w$ by a possible factor of $-1$, and
\begin{equation*}
\tau_w\tau_{w'}=\begin{cases}\pm\tau_{ww'}&\ell(w)+\ell(w')=\ell(ww'),\\0&\text{otherwise.}\end{cases}
\end{equation*}
\item The algebra $\ONH_n$ is a free $\Bbbk$-module.  Either of the sets $\lbrace\tau_wx^\alpha:w\in \mf{S}_n,\alpha_i\in\Z_{\geq0}\rbrace$, $\lbrace x^\alpha\tau_w:w\in \mf{S}_n,\alpha_i\in\Z_{\geq0}\rbrace$ is a basis.

\item 
 $\grdim \ONH_n = q^{-n(n-1)/2} [n]! \frac{1}{(1-q^2)^n}$
\end{enumerate}\end{proposition}

\subsubsection{Skew polynomial action}



Let $\Bbbk$ be a commutative ring (usually we take $\Bbbk=\Z$ or a field) and let
\begin{equation*}
\OPol_n=\Bbbk\langle x_1,\ldots,x_n\rangle/(x_ix_j+x_jx_i\text{ if }i\neq j)
\end{equation*}
be the $\Z$-graded superalgebra of \emph{skew polynomials} in $n$ variables.  Its generators are given degree $\deg(x_i) =2$, and parity $|x_i|=1$.  In particular, $$\grdim(\OPol_n) =  \left(\dfrac{1}{1-q^{2}}\right)^n.$$     The symmetric group $\mf{S}_n$ acts on $\OPol_n$ by
\begin{equation} \label{eq:sym-act}
w(x_i)=(-1)^{\ell(w)}x_{w(i)},\quad w(fg)=w(f)w(g)
\end{equation}
where $\ell(w)$ denotes the length of $w\in \mf{S}_n$. For $f \in \OPol_n$ we sometimes write $f^w$ for the action of $w$ on $f$.
For $i=1,\dots,n-1$, 
define the $i$-th \emph{odd divided difference operator} $\del_i$ to be the map $\OPol_n\to\OPol_n$ defined by
\begin{eqnarray} \label{eq:div-difference}
&\del_i(x_j)=\begin{cases}1&j=i,i+1\\0&\text{otherwise,}\end{cases}
\qquad
\del_i(fg)=\del_i(f)g+ s_i(f)\del_i(g).
\end{eqnarray}
For any $f\in \OPol_n$ the action of the odd divided difference operator is given by the formula
\[
 \partial_i(f) = \frac{(x_{i+1} - x_i)f -  s_i(f) (x_{i+1} - x_i) }{x_{i+1}^2-x_i^2},
\]
see \cite[equation 4.19]{KKO}.
These operators play a role analogous to that of the divided difference operators of Kostant-Kumar.

\begin{proposition}[Proposition 2.11 \cite{EKL} ]
There is a faithful action of $\ONH_n$ on $\OPol_n$ with $x_i$ acting by multiplication and $\tau_i$ acting by odd divided difference operators $\partial_i$.
\end{proposition}

\subsection{Odd symmetric functions}

Define the ring of \textit{odd symmetric polynomials} to be the subring
\begin{equation}
\OL_n=\bigcap_{i=1}^{n-1}{\rm Ker}(\partial_i)=\bigcap_{i=1}^{n-1}{\rm Im}(\partial_i)
\end{equation}
of $\OPol_n$.  In the even case, with the usual divided difference operators, this definition would agree with the usual notion of symmetric polynomials $\Lambda_n = \Z[x_1,\dots , x_n]^{\mf{S}_n}$.
Define elements of $\OPol_n$ for each $k\geq 1$
\begin{align}\label{eqn-defn-e}
\varepsilon_k(x_1,\ldots,x_n) &=\sum_{1\leq i_1<\cdots<i_k\leq n}\widetilde{x}_{i_1}\cdots\widetilde{x}_{i_k},\qquad\text{where }\widetilde{x}_i=(-1)^{i-1}x_i,
\\
h_k(x_1,\ldots,x_n) &=\sum_{1\leq i_1\leq\cdots\leq i_k\leq n}\widetilde{x}_{i_1}\cdots\widetilde{x}_{i_k},\qquad\text{where }\widetilde{x}_i=(-1)^{i-1}x_i,
\end{align}
Set $\varepsilon_0 = h_0 =1$ and $\varepsilon_j = h_j = 0$ for $j<0$.  If $j > n$ then $\varepsilon_j=0$.  Both of these families of skew polynomials are odd symmetric.  We call the $\varepsilon_i$ the $i$th \emph{odd elementary symmetric polynomial} and the $h_i$ the $i$th \emph{odd complete symmetric polynomial}.  The odd elementary and odd complete symmetric functions are related by
\begin{equation} \label{eq:odd-elemcomp}
  \sum_{j=0}^{\ell} (-1)^{\frac{1}{2}j(j+1)} \varepsilon_j h_{\ell-j} =0\quad \text{if $\ell \geq 1$}.
\end{equation}

It was shown in ~\cite{EKL} that the following relations hold in the ring $\OL_n$:
\begin{equation}\label{eqn-e-relations}\begin{split}
&\varepsilon_i\varepsilon_{2m-i}=\varepsilon_{2m-i}\varepsilon_i\qquad(1\leq i,2m-i\leq n)\\
&\varepsilon_i\varepsilon_{2m+1-i}+(-1)^i\varepsilon_{2m+1-i}\varepsilon_i
=(-1)^i\varepsilon_{i+1}\varepsilon_{2m-i}+\varepsilon_{2m-i}\varepsilon_{i+1}\qquad(1\leq i,2m-i\leq n-1)\\
&\varepsilon_1\varepsilon_{2m}+\varepsilon_{2m}\varepsilon_1=2\varepsilon_{2m+1}\qquad(1<2m\leq n-1).
\end{split}\end{equation}
Note that the third is the $i=0$ case of the second.  In particular, the ring $\OL_n$ of odd symmetric functions is noncommutative.

  For a partition $\l = (\l_1,\dots , \l_n)$ with $\l_1 \geq \l_2 \geq \dots \geq \l_n$ define the \emph{odd Schur polynomial}
\[
s_{\lambda}(x_1, \dots, x_n) := (-1)^{\binom{n}{3}} \partial_{w_0}(x^{\lambda} x_n^{\delta} x^\delta)^{w_0} = (-1)^{\binom{n}{3}} \partial_{w_0}(x_1^{\lambda_1} \dots x_n^{\lambda_n} x_1^{n-1} \dots x_n^{0} )^{w_0}
\]
where the superscript $w_0$ is the action of the longest element of $\mathfrak{S}_n$ defined in \eqref{eq:sym-act}.
Odd Schur polynomials give a $\Z$ basis of odd symmetric functions~\cite{EK,EKL}.  There is also an odd analog of the Littlewood-Richardson rule
\begin{equation} \label{eq:LR-Schur-mult}
 s_{\mu} s_{\nu} = \sum_{\l} c_{\mu \nu}^{\lambda} s_{\lambda}
\end{equation}for odd Schur polynomials developed by Ellis~\cite{Ellis-LR},
where the mod 2 reduction of $c_{\mu \nu}^{\lambda}$ agrees with the mod 2 reduction of the usual Littlewood-Richardson coefficients.  In particular, if $|\mu| + |\nu| \neq |\lambda|$, then $c_{\mu \nu}^{\lambda}=0$.

\begin{proposition}{\cite[Section 2.1]{EKL}}
The superalgebra $\OL_n$ has a presentation by generators $\varepsilon_1, \dots, \varepsilon_n$ and relations \eqref{eqn-e-relations}.  A basis of $\OL_n$ in $\Z$-degree $2j$ is given by all products
\begin{equation}
  \varepsilon_{\l} := \varepsilon_{\lambda_1} \dots \varepsilon_{\lambda_r}
\end{equation}
with $n\geq \lambda_1 \geq \lambda_2 \geq \dots \lambda_r \geq 1$ and $\lambda_1 + \dots + \lambda_r$ =j.  The same result holds if all $\varepsilon$'s are replaced by complete symmetric polynomials $h$'s.
\end{proposition}
%
%

Combining the results of \cite{EK} and the identification of various definitions of odd Schur functions from~\cite{Ellis-LR}, it follows that
\begin{align}
\label{eq:schur-complete}
  s_{\lambda} = h_{\lambda} + \sum_{\mu > \lambda} a_{\mu} s_{\mu}
   = \varepsilon_{\bar{\lambda}} + \sum_{\mu > \lambda} b_{\mu} s_{\mu}
\end{align}
for integers $a_{\mu}$ and $b_{\mu}$, where the order is the usual lexicographic order on partitions and $\bar{\lambda}$ denotes the dual partition.  We write $\ell(\l)$ for the number of nonzero parts of $\l$.

\begin{corollary}
$\grdim(\OL_m) = \prod_{i=1}^m \dfrac{1}{1-q^{2i}}$
\end{corollary}

\subsection{Odd Schubert polynomials}

For $w = s_{i_1} \dots s_{i_m}$ in $\mathfrak{S}_n$  define $\partial_w= \partial_{i_1} \dots \partial_{i_m}$.
 For the longest element $w_0$ of $\mathfrak{S}_n$ we fix a particular choice of reduced expression,
\begin{equation*}
\partial_{w_0}=\partial_1(\partial_2\partial_1)\cdots(\partial_{n-1}\cdots\partial_1).
\end{equation*}
For each $w\in \mathfrak{S}_n$ define the corresponding \textit{odd Schubert polynomial} $\sch_w\in\OPol_n$ by
\begin{equation} \label{eq-schubert}
\sch_w(x_1,\ldots,x_n)=\partial_{w^{-1}w_0}(\underline{x}^{\delta_n}).
\end{equation}
 For $w,w'\in \mathfrak{S}_n$, the formula
\begin{equation}\label{eqn_compose_oddops}
\partial_w\partial_{w'}=\begin{cases}\pm\partial_{ww'}&\text{if }\ell(ww')=\ell(w)+\ell(w')\\
0&\text{otherwise,}\end{cases}
\end{equation}
implies the following action of odd divided difference operators on odd Schubert polynomials
\begin{equation}\label{eqn_oddop_schubert}
\partial_u\sch_w=\begin{cases}\pm \sch_{wu^{-1}}&\text{if }\ell(wu^{-1})=\ell(w)-\ell(u)\\
0&\text{otherwise.}
\end{cases}\end{equation}
Note that the exact signs in these relations can be determined for specific choices of reduced words for the symmetric group elements involved, but we will not need these coefficients in what follows.

\begin{proposition}[Proposition 2.13 \cite{EKL}] \label{prop:OPol-basis}
$\OPol_n$ is a free left and right $\OL_n$-supermodule of graded rank  $q^{n(n-1)/2}[n]_{q,\pi}!$  with homogeneous graded bases give by either of the following sets
\begin{itemize}
  \item $\{ x_1^{a_1}x_2^{a_2}\dots x_n^{a_n} \; \mid \; a_i \leq n-i \;\; \text{for $1 \leq i \leq n$} \}$
  \item $\{\sch_w(x_1,\ldots,x_n) \; \mid \; w\in \mathfrak{S}_n \;   \}$
\end{itemize}
\end{proposition}

By \cite[Corollary 2.14]{EKL}this action on of the odd nilHecke algebra on odd polynomials   defines a superalgebra isomorphism
\begin{equation} \label{eq:varphi}
 \phi \maps \ONH_n \rightarrow \End_{\OL_n}(\OPol_n).
\end{equation}

\subsection{Modular reduction }


In what follows we make use of results from \cite[Section 2]{EKL} about various reductions mod 2.   Consider the reduction map $\symL{n} \to \symL{n} \otimes_{\Z} \Z_2$.  The images  of the usual elementary symmetric functions   under this map are nonzero in $\symL{n} \otimes_{\Z} \Z_2$.   In particular, $
\qrk(\symL{n})=\mathrm{grdim} (\symL{n} \otimes_{\Z} \Z_2)$.

Likewise, in the odd setting we also have
\begin{equation} \nn
\qrk(\OL_n)=\mathrm{grdim} (\OL_n \otimes_{\Z} \Z_2) \quad \text{and} \quad \qrk(\OPol_n)=\mathrm{grdim} (\OPol_n \otimes_{\Z} \Z_2).
\end{equation}
Products of odd elementary symmetric functions provide a $\Z_2$-basis for $\OL_n \otimes_{\Z} \Z_2$.
Since the definitions of odd divided difference operators, odd elementary symmetric functions, and odd polynomials agree with their even counterparts, when reduced modulo 2, we have isomorphisms
\begin{equation} \label{eq:iso_red_mod2}
\OPol_n\otimes_\Z \Z_2 \cong\Pol{n}\otimes_\Z \Z_2  \quad \text{and} \quad
\OL_n\otimes_\Z \Z_2 \cong\symL{n}\otimes_\Z \Z_2 .
\end{equation}
The first one fixes the generators $x_1,\ldots,x_n$, and the second one is the restriction of the first.


\begin{lemma}\label{lem:red_mod2}
	Let $A$ be an ring and let $M$ be an $A$-module such that $A, M$ are free as a $\Z$-modules. Let $S$ be a subset of $M$. Assume that $S\otimes 1$ is a linearly independent subset of $M\otimes_{\Z} \Z_2$ as an $A\otimes_{\Z} \Z_2$-module. Then $S$ is linearly independent over $A$.
\end{lemma}

\begin{proof}
	We proceed by contradiction. Assume that we have a non-trivial relation
	\begin{equation}\label{eq:dependence}
		\sum_{s \in S} a_ss=0.
	\end{equation}
	If all the coefficients in (\ref{eq:dependence}) are in $2A$, then we can divide the relation by $2$ to produce a new non-trivial relation by freeness of $M$. By freeness of $A$, this procedure terminates and we may assume that at least one $a_s$ is not in $2A$. In $M\otimes_{\Z} \Z_2$, we obtain a relation
	\[
		\sum_{s \in S} (a_s\otimes 1)(s\otimes 1)=0
	\]
	with at least one coefficient $a_s\otimes 1$ non-zero, a contradiction.
\end{proof}

\begin{corollary}\label{cor:inj_mod2}
	Let $A$ be an algebra free over $\Z$ and let $M,N$ be free $A$-modules. Let $f:M\rightarrow N$ be an $A$-linear map. Assume that $f\otimes_{\Z} \Z_2$ is injective. Then $f$ is injective.
\end{corollary}

\begin{proof}
	Let $S$ be a basis of $M$ over $A$. It suffices to prove that $f(S)$ is a linearly independent subset of $N$ over $A$. Since $f\otimes_{\Z} \Z_2$ is an injective $A\otimes_{\Z}\Z_2$-linear map and $S\otimes 1$ is a basis of $M\otimes_{\Z} \Z_2$ over $A\otimes_{\Z} \Z_2$, $\left(f\otimes_{\Z} \Z_2\right)(S\otimes 1)=f(S)\otimes 1$ is a linearly independent subset of $N\otimes_{\Z} \Z_2$ over $A\otimes_{\Z} \Z_2$. By Lemma \ref{lem:red_mod2}, it follows that $f(S)$ is a linearly independent subset of $N$ over $A$.
\end{proof}


\subsection{Deformed cyclotomic quotient}
Let $A_n$ be the graded supercommutative superalgebra generated by $\chi_1,\ldots,\chi_{\flr{\frac{n}{2}}}$  of parity $\bar{0}$ and $d$ of parity $\bar{1}$, subject to the relations
\begin{eqnarray}
	& & d^2 = 0, \qquad
	 d\chi_{\flr{\frac{n}{2}}} = 0 \quad \text{if } n \text{ is even}.
\end{eqnarray}
The $\Z$-grading on $A_n$ is defined by declaring $d$ to have degree 1 and $\chi_i$ to have degree $2i$. Note that for all odd elements $x,y$ of $A_n$, $xy=0$. Consequently, $A_n$ is also commutative. Put
\[
	c_i = \left \{ \begin{array}{ll}
		\chi_{\frac{i}{2}} & \text{ if } i \text{ is even},\\
		d\chi_{\frac{i-1}{2}} & \text{ if } i \text{ is odd}.
	\end{array} \right.
\]
The element $c_i$ is of degree $i$ and parity $\bar{i}$. Furthermore, we have the relation $c_{2i+1} = c_1c_{2i}=c_{2i}c_1$. \\
We can realize $A_n$ as a quotient of $\OL_n$, via the morphism
\begin{equation}\label{eq:an_quotient_map}
	\funct{\OL_n}{A_n,}{\varepsilon_i}{(-1)^{\binom{i}{2}}c_i.}
\end{equation}
 We will denote by $\bar{f}$ the image of an element $f\in \OL_n$ under this morphism.

 Let $A_n[t]:=A_n\otimes_\Bbbk \Bbbk[t]$ be the tensor product of superspaces where $t$ is parity $\bar{1}$.  For notational simplicity, in what follows we often omit tensor product for convenience an write $a f(t)$ in place of $a \otimes f(t)$ for $a\in A_n$ and $f(t) \in \Bbbk[t]$.  For instance, we write $a$ for $a\otimes 1$ and $f$ for $1 \otimes f$ for $a\in A$ and $f\in \Bbbk[t]$.

Consider the polynomial $a^n(t) \in A_n[t]$ defined by
\begin{equation} \label{eq:an}
	a^n(t)  = \sum_{\ell = 0}^{n} (-1)^{\ell} t^{\ell}c_{n-\ell}.
\end{equation}
An important point is that since $t$ has parity $\bar{1}$, the polynomial $a^n(t)$ neither commutes nor supercommutes with $t$. Instead, we have the relations

\begin{equation}\label{eq:an_relation}\begin{split}
		&ta^n(t) = a^n(t)(t-2d),\\
		&a^n(t)t = (t + (-1)^n2d)a^n(t).
\end{split}\end{equation}

\begin{definition}
The \emph{deformed cyclotomic quotient} is the graded superalgebra $\ONH_k^n$ defined as
\[
	\ONH_k^n = (\ONH_k \otimes A_n) / a^n(x_1).
\]
\end{definition}
The deformed cyclotomic quotients of Kang-Kashiwara-Oh \cite{KKO,KKO2} are specializations of $\ONH_k^n$. More precisely, we have an isomorphism of graded superalgebras
\[
	\ONH_k^n \otimes_{A_n} (A_n/d) \simeq R^n(k),
\]
where $R^n(k)$ is the notation used in \cite{KKO}.

\subsection{Odd equivariant cohomology of Grassmannians }

Recall that in the even case, the deformed cyclotomic quotient $\NH_k^n$ is isomorphic to a matrix ring over the $GL(n)$-equivariant cohomology ring $H^*_{GL(n)}(Gr(k;n))$ over the complex Grassmannian $Gr(k;n)$ of $k$-planes in $\C^n$,
\begin{equation}
  \NH_k^n \cong {\rm Mat}(k!, H^*_{GL(n)}(Gr(k;n)))
\end{equation}
see \cite[Proposition 9]{RoseWed}.  There is an isomorphism
\[
H^*_{GL(n)}(Gr(k;n)) \cong \Lambda_k \otimes \Lambda_{n-k}
\]
making $H^*_{GL(n)}(Gr(k;n))$ a free graded module of rank   $q^{k(n-k)}\qbin{n}{k}$   over $\Lambda_n$, with an integral basis given by the set of Schur polynomials $s^{\rm ev}_{\lambda} \in \Lambda_k$ with $\lambda=(\lambda_1,\dots, \lambda_k)$ and $\lambda_i \leq (n-k)$ for $1 \leq i \leq k$. In these statements, we view $\Lambda_k \otimes \Lambda_{n-k}$ as a subalgebra of $\Pol{n}$ via \eqref{eq:OLkOLnk2Opol}.

By analogy with the even case that will become apparent in this section, we
 define the \emph{odd equivariant Grassmannian bimodule}  as the $(\OL_k,A_n)$-bimodule
\begin{equation}
 OH^*_{GL(n)}(Gr(k;n))   := 
 (\OL_k \otimes \OL_{n-k}) \otimes_{\OL_n} A_n
\end{equation}
where we twist the right action of $\OL_n$ on $\OL_k \otimes \OL_{n-k}$ so that
\begin{equation}\label{eq:new-right-action}
  (a\otimes b) \cdot \varepsilon_r :=  (a\otimes b) \left( \sum\limits_{j=0}^r (-1)^{(r-j)j} \varepsilon_{r-j}\otimes \varepsilon_j\right).
\end{equation}

%
%

\begin{lemma} \label{lem:odd-eqcoh-basis}
The odd equivariant bimodule
$OH^*_{GL(n)}(Gr(k;n))$ is a free graded   $A_n$-supermodule of rank  $q^{k(n-k)}\qbin{n}{k}_{q,\pi}$   with basis given by the set
\begin{equation}
\cal{S} =\left\{
s_{\l} \otimes 1 \otimes 1 \mid \; \lambda=(\lambda_1,\dots, \lambda_k), \ \lambda_i \leq (n-k) \text{ for } 1 \leq i \leq k
\right\} .
\end{equation}
Alternatively, it has a basis given by
\begin{equation}
\cal{S}' = \left\{
1 \otimes s_{\mu} \otimes 1 \mid \; \lambda=(\lambda_1,\dots, \lambda_{n-k}), \ \lambda_i \leq (k) \text{ for } 1 \leq i \leq n-k
\right\}
\end{equation}
\end{lemma}

\begin{proof}
By base change, it suffices to prove that the set $\cal{S}$, respectively $\cal{S}'$,
is a basis of $\OL_k\otimes \OL_{n-k}$ as a right $\OL_n$-module.

The odd Schur polynomial  $s_\lambda$ acts on $\OL_k \otimes \OL_{n-k}$ via \eqref{eq:new-right-action}.  By  \eqref{eq:schur-complete} and     \eqref{eq:LR-Schur-mult} this is given by right multiplication
\begin{equation}
	(a\otimes b) \cdot s_{\lambda} = (a\otimes b) \cdot \sum_{\mu, \nu} \alpha_{\mu \nu}^{\lambda} s_{\mu}\otimes s_{\nu}
\end{equation}
for some coefficients $\alpha_{\mu \nu}^{\lambda} \in \Z$   equal to the odd Littlewood-Richardson coefficients mod 2. Hence, arguing inductively  as in the even case, we deduce that $\cal S$ spans $\OL_k\otimes \OL_{n-k}$ as a right $\OL_n$-module.

 We now prove that $\cal S$ is linearly independent over $\OL_n$. The isomorphisms \eqref{eq:iso_red_mod2} yield an isomorphism of algebras $\OL_n\otimes_{\Z}\Z_2\simeq \Lambda_n\otimes_{\Z} \Z_2$ and an isomorphism of right $\Lambda_n\otimes_{\Z} \Z_2$-modules
\begin{equation}\label{eq:iso_mod2}
\left(\OL_k\otimes \OL_{n-k}\right)\otimes_{\Z} \Z_2 \simeq \left(\Lambda_k\otimes \Lambda_{n-k}\right)\otimes_{\Z} \Z_2.
\end{equation}
Denote by $\cal S^{\rm ev}$ the even counterpart of $\cal S$. Since even and odd Schur functions coincide modulo 2, the isomorphism in (\ref{eq:iso_mod2}) sends $\cal S\otimes 1$ to $\cal S^{\rm ev}\otimes 1$. Furthermore, $\cal S^{\rm ev}$ is a basis of $\Lambda_{k}\otimes\Lambda_{n-k}$ as a right $\Lambda_n$-module, so $\cal S^{\rm ev}\otimes 1$ forms a basis of $\left(\Lambda_k\otimes \Lambda_{n-k}\right)\otimes_{\Z} \Z_2$ as a right $\Lambda_n\otimes_{\Z} \Z_2$-module. It follows that $\cal S\otimes 1$ is a basis of $\left(\OL_k\otimes \OL_{n-k}\right)\otimes_{\Z} \Z_2$ as a right $\OL_n\otimes \Z_2$-module. The conclusion now follows from Lemma \ref{lem:red_mod2}.
\end{proof}

\subsection{Deformed quotient as a matrix ring over the odd cohomology of a Grassmannian}

In this section we give a characterization of the deformed cyclotomic quotient as matrix ring over the odd equivariant cohomology of a Grassmannian. This description will be useful to study the braiding complexes in Section \ref{sec:braid-compl}.
%

Let 
$\phi \maps \ONH_k \rightarrow \End_{\OL_k}(\OPol_k)$
be the isomorphism \eqref{eq:varphi}. 
Under this isomorphism, a polynomial $f$ acts by multiplication.     To compute a set of defining relations for the odd cyclotomic quotient we must compute the image of the polynomial $a^n(x_1)$ from \eqref{eq:an}, and in particular, the matrix representing various powers of $\phi(x_1)$.  It is most convenient to work in the basis
\[
\cal{H}_{\alpha} := \{ \xt^{\alpha} \mid \alpha \leq (k-1,k-2, \dots, 1,0) \}
\]
of $\OPol_k$ as a free module over $\OL_k$. For each multi-index $\beta$ obtained by replacing $\alpha_1$ by zero in some $\alpha$ appearing in $\cal{H}_{\alpha}$, consider the $\OL_k$-submodule of $\OPol_k$ with basis
\[
 B_{\beta} = \{ \xt_1^{k-1}\xt^{\beta}, \xt_1^{k-2} \xt^{\beta}, \dots, \xt^{\beta} \}
\]
with $\xt^{\beta} = x_2^{\alpha_2} x_3^{\alpha_3} \dots x_k^{\alpha_k}$.
Then $\phi(x_1)$ sends the span of each $B_{\beta}$ to itself and is given in this basis by the
$k\times k$-matrices $X$    \cite[Lemma 5.1]{EKL},
where
\[X \; :=\; \phi(x_1)_\beta=
\begin{pmatrix}
\varepsilon_1 & 1 & 0 & 0 & \dots & 0\\
\varepsilon_2 & 0 & 1 & 0 & \dots & 0\\
-\varepsilon_3 & 0 & 0 & 1 & \dots & 0\\
\vdots & \vdots & \vdots & \vdots & \vdots & \vdots\\
(-1)^{\binom{k-2}{2}}\varepsilon_{k-1} & 0 & 0 & 0 & \dots & 1\\
(-1)^{\binom{k-1}{2}}\varepsilon_{k} & 0 & 0 & 0 & \dots & 0
\end{pmatrix}\]

\begin{lemma}
  Denote the $(i,j)$ entry of $X^m$ by $b_{i,j}^m$. Then $b_{i,j}^m$ are completely determined by the relations
  \begin{equation}\label{fmla:fruitloops}
    b_{i,j}^m=\left\{
\begin{array}{cl}
  h_{m+i-j}-\sum\limits_{l=1}^{i-1} h_l b_{i-l,j}^m  \quad &\text{if $j\leq m$, } \\
  \delta_{i+m,j}  \quad &\text{if $j> m$.} \\
\end{array}
\right.
  \end{equation}
\end{lemma}

\begin{proof}
We prove this by induction on $m$.
For $m=1$, the only non-trivial cases we need to check are for $j=1$ which follows from \eqref{eq:odd-elemcomp}
  \begin{align*}
		b_{i,1}^1 := (-1)^{\binom{i-1}{2}}\varepsilon_i
=h_{i}+\sum\limits_{l=1}^{i-1}(-1)^{\binom{l+1}{2}}h_{i-l}\varepsilon_l
=h_{1+i-1}-\sum_{l=1}^{i-1}(-1)^{\binom{i-l-1}{2}}h_l \varepsilon_{i-l}
	\end{align*}
For the inductive step,
observe by the exact same reasoning as \cite[Lemma 10]{RoseWed},   that $b_{i,j+1}^{m+1}=b_{i,j}^m$.
This implies that $b_{i,j}^{m+1}$ satisfy the required relations for $j\geq 2$, so all that remains to be checked are the $j=1$ entries.

We prove that the $b_{i,1}^{m+1}$ entries satisfy \eqref{fmla:fruitloops} by induction on $i$.
The induction hypothesis  on $m$ implies $b_{1,j}^{m}=h_{m+1-j}$. Thus, the $(1,1)$ entry in $X^{m+1}$ is
\[
b_{1,1}^{m+1}=\sum\limits_{j=1}^{k}(-1)^{\binom{j-1}{2}}h_{m+1-j}\varepsilon_j
=\sum\limits_{j=1}^{m}(-1)^{\binom{j-1}{2}}h_{m+1-j}\varepsilon_j=h_{m+1}.
\]
This proves the $i=1$ base case.
For the $i$-induction step, assume $b_{p,j}^{m+1}$ satisfies the relations for all $p<i+1$.
Then
\begin{align*}
  b_{i+1,1}^{m+1}&=\sum\limits_{r=1}^{k}(-1)^{\binom{r-1}{2}}b_{i+1,r}^{m}\varepsilon_r
  =\sum\limits_{r=1}^{k}(-1)^{\binom{r-1}{2}}\left(h_{m+i+1-r}-\sum\limits_{l=1}^{i}h_l b_{i+l-l,r}^m\right)\varepsilon_r
\\
  &=h_{m+(i+1)}-\sum\limits_{l=1}^{i}h_l\sum\limits_{r=1}^{k}(-1)^{\binom{r-1}{2}}b_{i+1-l,r}^m\varepsilon_r
  =h_{m+1+(i+1)-1}-\sum\limits_{l=1}^{i}h_lb_{i+1-l,r}^{m+1}
\end{align*}
Thus, $b_{i+1,j}^{m+1}$ satisfies the relations.
\end{proof}

The deformed cyclotomic quotient has defining relations in $\End_{\OL_k}(\OPol_k)$ given by setting the following element to zero
\[
\phi(a^n(x_1))=a^n(\phi(x_1))=a^n(X)=\sum\limits_{r=0}^n(-1)^{n-r} X^{n-r}\otimes c_r \in {\rm Mat}_{k!}(\OL_k \otimes A_n) .
\]
Denoting this matrix by $C$, we have that $C_{i,j}=(-1)^n\sum\limits_{r=0}^n (-1)^r b_{i,j}^{n-r}\otimes c_r$.
Then, under the isomorphism $\phi$ from \eqref{eq:varphi}
\begin{equation} \label{eq:deformed1}
\ONH_k^n \cong {\rm Mat}_{k!}
\left(\OL_k \otimes A_n / \langle C_{i,j} \mid 1\leq i,j \leq k \rangle \right).
\end{equation}
In particular, $(-1)^n C_{1,j}=\sum\limits_{r=0}^n (-1)^r h_{n+1-j-r}\otimes c_r=h_{n-j+1}'\in \OL_k \otimes A_n$, where we define
\begin{equation} \label{eq:hm'}
  h_{m}' =\sum\limits_{r=0}^m (-1)^r h_{m-r}\otimes c_r \in \OL_k \otimes A_n .
\end{equation}

\begin{lemma} \label{lem:I'}
 Let $I'$ denote the ideal generated by the entries of the first row of $C$,
\begin{align*}
  I':=\langle C_{1,j} \mid 1\leq j \leq k \rangle = \langle h_m'  \mid n-k+1 \leq m \leq n \rangle.
\end{align*}
Then   $C_{i,j}\in I'$ for all $1\leq i,j \leq k$.
\end{lemma}

\begin{proof}
  We prove this by induction on $i$. The statement is immediate by definition for $i=1$, so now assume that $C_{p,j}\in I'$ for all $1\leq p <i$.
  Then we have
  \begin{align*}
    (-1)^nC_{i,j}&=\sum\limits_{r=0}^{n}(-1)^{r}\left(h_{n+1-r-j} -\sum\limits_{l=1}^{i-1} h_l b_{i-l,j}^{n-r} \right)\otimes c_r =h_{n+i-j}' -\sum\limits_{l=1}^{i-1} h_l C_{i-l,j}
  \end{align*}
  Each $ h_l C_{i-l,j}$ is in $I'$ by the induction hypothesis, so we only need to show that $h_{n+i-j}' \in I'$ for all $1\leq i,j\leq k$.
  We know that $h_{n+i-j}'$ is a generator of $I'$ for $i< j$, so it suffices to prove that $h_{n+s}' \in I'$ for all $s\geq 0$.
  We prove this with a simple induction argument on $s$.
  The case $s=0$ is clearly true, so now assume that $h_{n+r}' \in I'$ for all $0\leq r < s $.
  It follows that
  \begin{align*}
    h_{n+s}'  &=\sum\limits_{r\geq 0}(-1)^r h_{n+s-r}\otimes c_r = -\sum\limits_{r\geq 0}(-1)^r\left(\sum\limits_{i\geq 1} (-1)^{\binom{i+1}{2}}\varepsilon_i h_{n+s-r-i}\otimes c_r\right) \\
&=-\sum\limits_{i\geq 1} (-1)^{\binom{i+1}{2}} \varepsilon_i  h_{n+s-i}'  \in I'
  \end{align*}
\end{proof}

Consider the $A_n$-superalgebra defined by
\begin{equation}
M_k^n :=  \left( \OL_k \otimes A_n \right)/\langle h_{m}' \mid n-k+1\leq m \leq n \rangle.
\end{equation}
For reasons that will become apparent, we call this the odd equivariant Grassmannian cohomology\footnote{We emphasize that `odd cohomology' is not an actual cohomology theory that can be applied to any manifold.  In general, `odd cohomology' is a noncommutative analog of usual singular cohomology that becomes isomorphic to the usual singular cohomology when coefficients are reduced modulo two.  See \cite{Lau-odd} for a discussion of odd cohomologies of Grassmannians and \cite{LauR} for a discussion on odd cohomologies of Springer varieties and a broader `oddification' program. } superalgebra.

\begin{theorem} \label{thm:matrixMnk}
There is an isomorphism of graded superalgebras 
\begin{equation}
 \ONH_k^n \cong {\rm Mat}_{k!}( M_k^n).
\end{equation}
This implies that the isomorphism $\phi \maps \ONH_n \rightarrow \End_{\OL_n}(\OPol_n)$ induces an isomorphism
\begin{equation}
  \ONH_k^n \longrightarrow
  \End_{M_{k}^n}(\OPol_k \otimes_{\OL_k} M_k^n).
\end{equation}
\end{theorem}

\begin{proof}
 The deformed cyclotomic quotient is given by \eqref{eq:deformed1}.  By Lemma~\ref{lem:I'} the quotient is generated by the first row.  The result follows.
\end{proof}

\begin{proposition} \label{prop:Mkn-action}
The odd equivariant Grassmannian bimodule  $OH^*_{GL(n)}(Gr(k;n))$ is a   $(M_k^n,M_k^n)$-bimodule with left action given by
\[
(f\otimes a) \cdot (1 \otimes s_{\mu} \otimes 1) =  (-1)^{|a||\mu|}(f  \otimes s_{\mu} \otimes a)
\]
and right action
\[
(1 \otimes s_{\mu} \otimes 1) \cdot (f\otimes a)=  (-1)^{|f||\mu|}(f  \otimes s_{\mu} \otimes a)
\]
\end{proposition}

\begin{proof}
Lemma~\ref{lem:odd-eqcoh-basis} gives a bases for  $OH^*_{GL(n)}(Gr(k;n))$ as an $A_n$-supermodule.  To see that this map is well defined, it suffices to show that $h_m'$ acts by $0$ for $m>n-k$.
\begin{align*}
 h'_m \cdot (1 \otimes s_{\mu} \otimes 1) &= \sum_i (-1)^i(h_{m-i} \otimes c_i) \cdot (1 \otimes s_{\mu} \otimes 1)
\\  &=\sum_{i=0}^m (-1)^{i+i|\mu |} (h_{m-i}  \otimes s_{\mu} \otimes c_i)
  =\sum_{i=0}^m \sum_{b=0}^i (-1)^{i+b|\mu| +\binom{i-b}{2}+\binom{b}{2}} (h_{m-i}\varepsilon_{i-b}  \otimes s_{\mu}\varepsilon_b \otimes 1)
  \\ &= \sum_{b\geq 0}(-1)^{b|\mu|+\binom{b}{2}+b}\sum_{i'=0}^{m-b} (-1)^{i'+\binom{i'}{2}} (h_{m-b-i'}\varepsilon_{i'}  \otimes s_{\mu}\varepsilon_b \otimes 1) = 0
\end{align*}
The last equality is true for $m>n-k$ by \eqref{eq:odd-elemcomp}.
\end{proof}

We will show in Proposition~\ref{prop:eqcoh} that $M_k^n$ is isomorphic to $OH^*_{GL(n)}(Gr(k;n))$ as  $(M_k^n, M_k^n)$-bimodules.

\begin{lemma} \label{lem:composition}
Let $\mu=(\mu_1,\dots, \mu_k)$ be a composition satisfying $\mu_i \leq n-k$ for all $1 \leq i \leq k$.  Then
the product $h_{\mu} := h_{\mu_1} \dots h_{\mu_k}$  of odd complete symmetric functions in $M_k^n$ is in the span of the set
 \begin{equation} \label{eq:odd-eq-complete0}
\{
h_{\lambda} \otimes a \mid a\in A_n; \;\; \lambda = (\lambda_1, \dots, \lambda_{k}), \; \;  \text{with $0 \leq \l_k \leq \dots \l_1 \leq (n-k)$}
\} .
\end{equation}
\end{lemma}

\begin{proof}
We prove this by induction on $|\mu|$. For $|\mu|=0$, this is trivial as $h_0=1\otimes 1$.
Assume then that $h_\mu$ is in the span of the set \eqref{eq:odd-eq-complete0} for all compositions $\mu=(\mu_1, \dots , \mu_k)$ with $\mu_i\leq n-k$ for all $1\leq i \leq k$ and $|\mu|<r$.
We establish the
$|\mu|=r$ case.

Let $h_{\nu}$ satisfy the hypothesis of the Lemma with $|\nu|=r$.  We use \cite[Remark 2.4]{EKL} with $\varepsilon$ replaced by $h$  to sort a product of odd complete symmetric functions into a non-increasing order.
If $p<l$ and $p+l$ is odd, then
\begin{align}\label{eq:h-shuffle}
  h_ph_l= \left\{\begin{array}{c}
                   h_lh_p+2\sum\limits_{i=1}^p (-1)^{\binom{i}{2}} h_{l+i}h_{p-i} \quad \text{if $p$ even,} \\
                   -h_lh_p+2\sum\limits_{i=1}^p (-1)^{\binom{i-1}{2}} h_{l+i}h_{p-i} \quad \text{if $p$ odd.}
                 \end{array}\right.
\end{align}
Hence, $h_\nu$ can be written
  $h_\nu=\sum\limits_{\gamma }c_\gamma h_\gamma$
for integers $c_{\gamma}$ where $\gamma=(\gamma_1, \dots, \gamma_k)$ varies over the set of partitions of $|\nu|=r$ with $\ell(\gamma)\leq \ell(\nu)$.
It suffices to show that all of the $h_{\gamma}$ lie in the span \eqref{eq:odd-eq-complete0}.

For any $a\in A_n$ and composition $\alpha=(\alpha_1,\dots \alpha_k)$ with $\alpha_i>n-k$, the relations of $M_k^n$ imply
\begin{align*}
  h_\alpha\otimes a &=\sum\limits_{j=1}^{\alpha_i} (-1)^{\alpha_i+1+j+j(\alpha_{i+1}+ \dots + \alpha_{k})} (h_{ (\alpha_{1}, \dots, \alpha_{i-1} ,\alpha_{i} -j,\alpha_{i+1}, \dots, \alpha_{k})}\otimes c_ja).
\end{align*}
Therefore, we can write each $h_\gamma$ as a linear combination
$
  h_{\gamma}= \sum\limits_{\beta} b_\beta h_\beta \otimes d_\beta
$
for some  $d_\beta \in A_n$, and $b_\beta \in \Z$ where $\beta=(\beta_1, \dots \beta_k)$ varies over compositions with $\beta_i\leqslant n-k$ for all $1\leq i \leq k$ and $|\beta|<|\nu|=r$.
Each $h_\beta$ lies in the span of \eqref{eq:odd-eq-complete0} by the induction hypothesis, and the result follows.
\end{proof}

\begin{lemma} \label{lem:Mkn-basis}
The $A_n$-supermodule $M_k^n$ is 
 spanned by  the set
\begin{equation} \label{eq:odd-eq-complete}
\{
h_{\lambda} \mid \lambda = (\lambda_1, \dots, \lambda_{k}), \; \;  \text{with $0 \leq \l_k \leq \dots \l_1 \leq (n-k)$}
\} .
\end{equation}
\end{lemma}

\begin{proof}
 The superalgebra $\OL_k$ has a basis given by
$\{ s_{\lambda}  \mid  \lambda = (\lambda_1, \dots, \lambda_k) \}$. It follows from \eqref{eq:schur-complete} and $\varepsilon_{m} =0 \in \OL_k$ for $m>k$ that any Schur function $s_{\l}$ with $\ell(\l)>k$ vanishes in $\OL_k$.
We argue that any $h_{\mu}$ with $\ell(\mu) >k$ can be written as a linear combination of $h_{\mu'}$ with $\ell(\mu')\leq k$ in $\OL_k$.  To see this, use \eqref{eq:schur-complete} to write $h_{\mu}$ in terms of Schur functions
\begin{equation} \label{eq:BB1}
 h_{\mu}  = s_{\mu} - \sum_{\mu' > \mu} a_{\mu} s_{\mu'} .
\end{equation}
The $s_{\mu}$ term vanishes in $\OL_k$ and if $\mu'$ has $\ell(\mu')>k$  this term also vanishes in the summation.  Otherwise, since $\mu' > \mu$ it suffices to assume that all $\mu'$ appearing in the summation have $\ell(\mu')\leq k$.
Now we use \eqref{eq:schur-complete}
so that \eqref{eq:BB1} can be written
\[
 h_{\mu}  =   - \sum_{\mu' > \mu} a_{\mu} s_{\mu'} = - \sum_{\mu' > \mu} a_{\mu}h_{\mu'} - \sum_{\mu' > \mu}\sum_{\mu'' > \mu'} a_{\mu} c_{\mu'} h_{\mu''}.
\]
Since $\mu' < \mu''$ are both partitions of the same size and $\ell(\mu') \leq k$, it follows that all $\mu''$ appearing in the sum have $\ell(\mu'') \leq k$.

Therefore, to show that the set \eqref{eq:odd-eq-complete} spans, it suffices to show that for any $k$-part partition $\lambda = (\l_1, \dots, \l_k)$, we can write $h_{\l}$ as an $A_n$-linear combination of elements from \eqref{eq:odd-eq-complete}.
Observe that  in the quotient module $M_k^n$, any $h_j=h_j\otimes 1$ with $j> n-k$ can be written as an $A_n$-linear combination of $h_{\ell}$ with $\ell \leq n-k$.  However, because of the noncommutativity of the odd complete symmetric functions, this rewriting process may produce terms $h_{\mu} \otimes a \in \OL_k \otimes A_n$ where $\mu$ is a composition $\mu=(\mu_1, \dots, \mu_k)$ satisfying $\mu_i \leq n-k$ for all $1 \leq i \leq k$.  The claim then follows from Lemma~\ref{lem:composition}. Hence, the set \eqref{eq:odd-eq-complete} is a spanning set for $M_k^n$.
\end{proof}

\begin{proposition} \label{prop:eqcoh}
 There is an isomorphism of $(M_k^n,M_k^n)$-bimodules
\begin{align}
\beta \maps M_k^n
&\rightarrow
OH^*_{GL(n)}(Gr(k;n))\\
\nn
f \otimes a & \mapsto f \otimes 1 \otimes a.
\end{align}
Therefore, the spanning set from Lemma~\ref{lem:Mkn-basis} of the $A_n$-supermodule $M_k^n$ is a basis.
\end{proposition}

\begin{proof}
Recall the $(M_k^n,M_k^n)$-bimodule structure on $OH^*_{GL(n)}(Gr(k;n))$ from Proposition~\ref{prop:Mkn-action}.
 To see that this map is well-defined, observe that $\beta(f\otimes a) = (f\otimes a) \cdot (1\otimes 1 \otimes 1)$.
%
%
 Now observe that for a partition $\lambda$, the image of $h_{\lambda}\otimes 1$ is
	\[
		\left(s_{\lambda} -  \sum_{\mu > \lambda} a_{\mu} s_{\mu}\right) \otimes 1 \otimes 1
	\]
	where the coefficients $a_{\mu}$ are as in \eqref{eq:schur-complete}. Therefore, the image of the spanning set (\ref{eq:odd-eq-complete}) from Lemma \ref{lem:Mkn-basis} is a basis of $OH^*_{GL(n)}(Gr(k;n))$ as an $A_n$-module by Lemma \ref{lem:odd-eqcoh-basis}. The conclusion follows.
\end{proof}


\subsection{Isomorphisms for complementary odd Grassmannian cohomology superalgebras}
Just as there are isomorphisms $H^*_{GL(n)}(Gr(k;n)) \cong H^*_{GL(n)}(Gr(n-k;n))$ between equivariant cohomology rings for complementary Grassmannians in the usual setting, in this section we show that the same holds for odd equivariant cohomology superalgebras $M_k^n \cong M_{n-k}^n$.

We give an alternative basis for $M_k^n$ that is more convenient for specifying the isomorphism.  Recall the notation from \eqref{eq:hm'}.
\begin{lemma} \label{lem:hpverson}
  The $A_n$-supermodule $M_k^n$ has a basis given by
  \[
  \{h_\l' \mid \l=(\l_1,\dots \l_k)\; \text{with $0 \leq \l_k \leq \dots \l_1 \leq (n-k)$} \}.
  \]
\end{lemma}

\begin{proof}
The change of basis from $h'$ to $h$ is unipotent.
By Proposition~\ref{prop:eqcoh} it follows that the set of $h'$'s form a  basis.
\end{proof}

\begin{lemma}  \label{lem:primeoddsym}
The $A_n$-superalgebra $M_k^n$ has a presentation by generators $\{h_r' \mid r \geq 0\}$ 
modulo relations
\begin{equation}\label{eqn-h'-relations}\begin{split}
&h_i'h_{2m-i}'=h_{2m-i}'h_i'\qquad(1\leq i,2m-i\leq n)\\
&h'_ih'_{2m+1-i}+(-1)^ih'_{2m+1-i}h'_i
=(-1)^ih'_{i+1}h'_{2m-i}+h'_{2m-i}h'_{i+1}\qquad(1\leq i,2m-i\leq n-1)\\
&h'_1 h_{2m}+h'_{2m}h'_1=2h'_{2m+1}\qquad(1<2m\leq n-1) 
\end{split}\end{equation}
and $h'_r=0$ for $r > n-k$.
\end{lemma}

\begin{proof}
The set $\{h_r' \mid r \leq n-k\}$ generates $M_k^n$ as a superalgebra.
One can check by direct computation that these relations hold in $M_k^n$.  To see that there are no additional relations in $M_k^n$, observe that using only these relations, we can write any product $h_{r_1}' \dots h_{r_m}'$ as a linear combination of elements with decreasing order, as in \cite[Remark 2.4]{EKL}, to give a linear combination of basis  elements from Lemma~\ref{lem:hpverson}.

\end{proof}

\begin{proposition} \label{prop:Mkn2Mnmkn}
There is an isomorphism of $A_n$-superalgebras $\gamma_k^n \maps M_k^n \to M_{n-k}^n$ defined by
\[
\gamma_k^n(h_r')=(-1)^{(k+1)r+\binom{r}{2}}\varepsilon_r\otimes 1 
\]
\end{proposition}

\begin{proof}
The elements $\gamma_k^n(h_r')$ satisfy the relations in Lemma~\ref{lem:primeoddsym} by \eqref{eqn-e-relations}.  Hence, this is a well defined map of superalgebras.   The set of elements $\{ \varepsilon_r \mid 0 \leq r \leq n-k\}$ generate $M_{n-k}^n$ as a superalgebra by \eqref{eq:odd-elemcomp} and Lemma~\ref{lem:Mkn-basis}, so the map is a surjection.  It follows that  $\gamma_k^n$ is an isomorphisms since $M_k^n$ and $M_{n-k}^n$ have the same graded dimension by Lemma~\ref{lem:odd-eqcoh-basis}.
\end{proof}

%
\section{The odd 2-category and faithful universal quotients}\label{subsec:definitionodd2cat}
%

\subsection{Definition of the odd 2-category}


We recall here the definition of the rank one super Kac-Moody 2-category from \cite{BE2}.  This presentation greatly simplifies the presentation from \cite{Lau-odd} where it was referred to as the odd categorification of $\mf{sl}_2$.

\begin{definition} \label{def:oddU}
The odd $2$-supercategory $\oUcat=\oUcat(\mf{sl}_2)$ is the $2$-supercategory
consisting of
\begin{itemize}
\item objects $\l$ for $\l \in \Z$,
\item generating 1-morphisms $\cal{E}\1_{\lambda} \maps \lambda \to \lambda+2$ and $\cal{F}\1_{\l}\maps \lambda \to \lambda -2$ for each $\l \in \Z$,

\item $\Z\times \Z_2$-graded generating 2-morphisms $x\maps \cal{E}\1_{\l} \to \cal{E} \1_{\l}$  of  parity $\bar{1}$ and degree $2$,  $\tau \maps \cal{E}\cal{E} \1_{\l} \to \cal{E}\cal{E} \1_{\l}$ of  parity $\bar{1}$ and degree -2, and $\eta \maps \1_{\l} \to \cal{F}\cal{E}\1_{\l}$ and $\varepsilon \maps \cal{E}\cal{F}\1_{\l} \to \1_{\l}$ both of parity $\bar{0}$ and degree $1+\l$ and $1-\l$, respectively.  These generators are subject to the relations given below, which are most conveniently expressed in the string diagrammatics of supercategories from \cite{BE1}.
\begin{align}
x:=
\hackcenter{\begin{tikzpicture}[scale=0.8]
    \draw[thick, ->] (0,0) -- (0,1.25)
        node[pos=.5, shape=coordinate](DOT){};
    \filldraw  (DOT) circle (2.5pt);
    \node at (-.85,.85) {\tiny $\lambda +2$};
    \node at (.5,.85) {\tiny $\lambda$};
\end{tikzpicture}}
\qquad
\tau := \hackcenter{\begin{tikzpicture}[scale=0.8]
\draw[thick, ->] (0,0) .. controls (0,.5) and (.75,.5) .. (.75,1.0);
\draw[thick, ->] (.75,0) .. controls (.75,.5) and (0,.5) .. (0,1.0);
\node at (1.1,.55) {\tiny $\lambda$};
\end{tikzpicture}}
\qquad
  \eta:= \hackcenter{\begin{tikzpicture}[scale=0.8]
  		\draw[thick, <-] (.75,2) .. controls ++(0,-.75) and ++(0,-.75) .. (0,2);
  		\node at (.4,1.2) {\tiny $\lambda$};
  \end{tikzpicture}}
 \qquad
\varepsilon :=
\hackcenter{\begin{tikzpicture}[scale=0.8]
	\begin{scope}[yscale=-1]
    \draw[thick, <-] (.75,2) .. controls ++(0,-.75) and ++(0,-.75) .. (0,2);
    \node at (.4,1.2) {\tiny $\lambda$}; \end{scope}
\end{tikzpicture}}
\end{align}
\end{itemize}
The identity $2$-morphism of the $1$-morphism $\cal{\cal{E}} \onen$ is represented by an upward oriented line (likewise, the identity $2$-morphism of $\cal{F} \onen$ is represented by a downward oriented line).

Horizontal and vertical composites of the above diagrams are interpreted using the conventions for supercategories explained in \cite{BE1}.   The  2-supercategory structure implies that diagrams with odd parity skew commute.
The $2$-morphisms satisfy the following relations
 (see \cite{BE2} for more details).
\begin{enumerate}

\item     Odd nilHecke: The odd nilHecke relations from Definition~\ref{def-oddnil} are satisfied for upward oriented strands and any $\l\in \Z$.

\item    Right adjunction axioms: $\cal{F}\1_{\l+2}$  is a right dual of $\cal{E}\1_{\l}$ with unit $\eta$ and counit $\varepsilon$.

\item Odd $\mf{sl}_2$ isomorphisms:  we define a new 2-morphisms
\begin{align} \label{eq:crossl-gen-cyc}
\hackcenter{
\begin{tikzpicture}[scale=0.8]
    \draw[thick, ->] (0,0) .. controls (0,.5) and (.75,.5) .. (.75,1.0);
    \draw[thick, <-] (.75,0) .. controls (.75,.5) and (0,.5) .. (0,1.0);
    \node at (1.1,.65) {\tiny $\lambda$};
\end{tikzpicture}}
  := \;\;
\hackcenter{\begin{tikzpicture}[scale=0.7]
    \draw[thick, ->] (0,0) .. controls (0,.5) and (.75,.5) .. (.75,1.0);
    \draw[thick, ->] (.75,-.5) to (.75,0) .. controls (.75,.5) and (0,.5) .. (0,1.0) to (0,1.5);
    \draw[thick] (0,0) .. controls ++(0,-.4) and ++(0,-.4) .. (-.75,0) to (-.75,1.5);
    \draw[thick, ->] (.75,1.0) .. controls ++(0,.4) and ++(0,.4) .. (1.5,1.0) to (1.5,-.5);
    \node at (1.85,.55) { \tiny $\lambda$};
\end{tikzpicture}}
 \maps \cal{E}\cal{F}\1_{\l} \to \cal{F}\cal{E}\1_{\l}
\end{align}
of parity $\bar{1}$ and degree $0$.  Then there are (non homogeneous) isomorphisms
\begin{align} \label{prop_coveringrelsU1}
& \hackcenter{
\begin{tikzpicture}[scale=0.8]
    \draw[thick, ->] (0,0) .. controls (0,.5) and (.75,.5) .. (.75,1.0);
    \draw[thick, <-] (.75,0) .. controls (.75,.5) and (0,.5) .. (0,1.0);
    \node at (1.1,.65) {\tiny $\lambda$};
\end{tikzpicture}}\;\;   \;\; \bigoplus_{k=0}^{\lambda-1}
\hackcenter{\begin{tikzpicture}[scale=0.8]
  \draw[thick,  <-] (0.35,.2) .. controls ++(0.1,0.8) and ++(-0.1,0.8) .. (-0.35,.2)  node[pos=0.80]{\tikz \draw[fill=black] circle (0.4ex);};
\node at (-.4,.9) {$\scs k$};
\end{tikzpicture} }
: \mathcal{E} \mathcal{F}  \1_{\lambda} \rightarrow
\mathcal{F}\mathcal{E}  \1_{\lambda}
\bigoplus_{k=0}^{\lambda-1}    \1_{\lambda}
 &\text{for $\l \geq 0$ } \\
 \label{prop_coveringrelsU2}
&\hackcenter{
\begin{tikzpicture}[scale=0.8]
    \draw[thick, ->] (0,0) .. controls (0,.5) and (.75,.5) .. (.75,1.0);
    \draw[thick, <-] (.75,0) .. controls (.75,.5) and (0,.5) .. (0,1.0);
    \node at (1.1,.65) {\tiny $\lambda$};
\end{tikzpicture}}\;\; \bigoplus_{k=0}^{-\l-1}
\hackcenter{\begin{tikzpicture}[scale=0.8]
  \draw[thick, ->] (-0.37,.5) .. controls ++(-.1,-.8) and ++(.1,-.8) .. (0.35,.5)
      node[pos=0.85, shape=coordinate](Y){}
      node[pos=0.55, shape=coordinate](M){}
      node[pos=0.44, shape=coordinate](X){};
   \node at (Y){\tikz \draw[fill=black] circle (0.4ex);};
\node at (.55,.15) {$\scs k$};
\end{tikzpicture} }
:\mathcal{E}\mathcal{F}   \1_{\lambda} \bigoplus_{k=0}^{-\l-1}
\1_{\lambda} 
\rightarrow \mathcal{F} \mathcal{E} \1_{\lambda}
 &\text{for $\l \leq 0$}
\end{align}
in the 2-category $\oUcat$.
\end{enumerate}
\end{definition}

Ellis and Brundan argue that this compact definition implies the existence of 2-morphisms $\eta' \maps \1_{\l} \to \cal{E}\cal{F}\1_{\l}$ and $\varepsilon' \maps \cal{F}\cal{E}\1_{\l} \to \1_{\l}$ satisfying left adjunction axioms up to a sign.  We depict these maps as
\begin{align}
 \eta' :=
\hackcenter{\begin{tikzpicture}[scale=0.8]
    \draw[thick, ->] (.75,2) .. controls ++(0,-.75) and ++(0,-.75) .. (0,2);
    \node at (.4,1.2) {\tiny $\lambda$};
\end{tikzpicture}}  & &
\varepsilon' := \hackcenter{\begin{tikzpicture}[scale=0.8]
    \draw[thick, ->] (.75,-2) .. controls ++(0,.75) and ++(0,.75) .. (0,-2);
    \node at (.4,-1.2) {\tiny $\lambda$};
\end{tikzpicture}}
 \\
 \quad \text{degree } (1-\l, \overline{\l+1})  & &
 \quad \text{degree } (1+\l,\overline{\l+1}) \nn
\end{align}
where we have indicated a $Q$-grading and parity as an ordered tuple $(x,\bar{y})$. One can also deduce that $\ONH_a^{{\rm sop}}$ acts on $\cal{F}^a\onel$ for all $\l$ and $a \geq 0$ where the generators $\tau_i^{{\rm sop}}$ and $x_i^{{\rm sop}}$ act by
\[
x^{{\rm sop}} \mapsto
\hackcenter{\begin{tikzpicture}[scale=0.8]
    \draw[thick, <-] (0,0) -- (0,1.25)
        node[pos=.5, shape=coordinate](DOT){};
    \filldraw  (DOT) circle (2.5pt);
    \node at (-.85,.85) {\tiny $\lambda +2$};
    \node at (.5,.85) {\tiny $\lambda$};
\end{tikzpicture}} :=
\hackcenter{\begin{tikzpicture}[scale=0.8]
    \draw[thick, ->]  (0,.4) .. controls ++(0,.6) and ++(0,.6) .. (.75,.4) to (.75,-1);
    \draw[thick, <-](0,.4) to (0,-.4) .. controls ++(0,-.6) and ++(0,-.6) .. (-.75,-.4) to (-.75,1);
    \filldraw  (0,-.2) circle (2.5pt);
    \node at (1.3,.9) {\tiny $\lambda + 2$};
\end{tikzpicture}}
\qquad
\tau^{{\rm sop}} \mapsto
\hackcenter{
\begin{tikzpicture}[scale=0.8]
    \draw[thick, <-] (0,0) .. controls (0,.5) and (.75,.5) .. (.75,1.0);
    \draw[thick, <-] (.75,0) .. controls (.75,.5) and (0,.5) .. (0,1.0);
    \node at (1.1,.65) {\tiny $\lambda$};
\end{tikzpicture}}
\;\; := \;\;
\hackcenter{\begin{tikzpicture}[scale=0.7]
    \draw[thick, ->] (0,0) .. controls (0,.5) and (.75,.5) .. (.75,1.0);
    \draw[thick, ->] (.75,0) .. controls (.75,.5) and (0,.5) .. (0,1.0);
    \draw[thick] (0,0) .. controls ++(0,-.4) and ++(0,-.4) .. (-.75,0) to (-.75,2);
    \draw[thick] (.75,0) .. controls ++(0,-1.2) and ++(0,-1.2) .. (-1.5,0) to (-1.55,2);
    \draw[thick, ->] (.75,1.0) .. controls ++(0,.4) and ++(0,.4) .. (1.5,1.0) to (1.5,-1);
    \draw[thick, ->] (0,1.0) .. controls ++(0,1.2) and ++(0,1.2) .. (2.25,1.0) to (2.25,-1);
    \node at (-.35,.75) {\tiny  $\lambda$};
\end{tikzpicture}}
\]

\subsection{Bubble Relations} \label{sec:bubbles}
 The isomorphisms \eqref{prop_coveringrelsU1}  -- \eqref{prop_coveringrelsU2} imply a number of other diagrammatic relations, see \cite{Lau-odd,BE2}.
In particular, these isomorphisms imply important relations for `dotted bubbles' which are  various composites of the 2-morphisms $\eta$, $\eta'$, $\varepsilon$, $\varepsilon'$ and $x$ giving endomorphisms of $\1_{\l}$.    The relations imply that dotted bubbles of negative degree are zero and degree zero bubbles are equal to the identity, so that
for all $m \geq 0$
\begin{equation}
 \hackcenter{ \begin{tikzpicture} [scale=.65]
 \draw[ ]  (0,0) arc (180:360:0.5cm) [thick];
 \draw[<- ](1,0) arc (0:180:0.5cm) [thick];
\filldraw  [black] (.1,-.25) circle (2.5pt);
 \node at (-.2,-.5) {\tiny $m$};
 \node at (1.15,.8) {\tiny $\lambda  $};
\end{tikzpicture} } \;  = \delta_{m,\l-1}\Id_{\1_{\l}}\quad \text{if $m \leq \l -1$}, \qquad \quad
\;
\hackcenter{ \begin{tikzpicture} [scale=.65]
 \draw  (0,0) arc (180:360:0.5cm) [thick];
 \draw[->](1,0) arc (0:180:0.5cm) [thick];
\filldraw  [black] (.9,-.25) circle (2.5pt);
 \node at (1,-.5) {\tiny $m$};
 \node at (1.15,.8) {\tiny $\lambda $};
\end{tikzpicture} } \;  = \delta_{m,-\l-1}\Id_{\1_{\l}} \quad  \text{if $m \leq -\l -1$}.
\end{equation}
We will sometimes make use of the shorthand notation that highlights the degree of the bubble which are $2n$ in both cases below
\begin{equation}
   \vcenter{\hbox{$\posbubd{}{ n+\ast\quad }$}} :=
     \vcenter{\hbox{$\posbubd{}{\l-1+n\qquad} $}}
\qquad
   \vcenter{\hbox{$\negbubd{}{\; \;n+\ast }$}} :=
     \vcenter{\hbox{$\negbubd{}{\qquad \;-\l-1+n\;} $}}
\end{equation}

 The degree two bubble is given a special notation as in \eqref{eq:def-odd-bubble} and squares to zero by the superinterchange law.  We call this map the odd bubble and denote it by
    \begin{equation} \label{eq:def-odd-bubble}
\raisebox{-4mm}{$\oddbubble{ }$}
\; :=\; \left\{
    \begin{array}{cl}
       \posbubd{}{\l}
 & \text{if $\l \geq 0$} \\
        \negbubd{}{-\l}  & \text{if $\l\leq 0$}
    \end{array}
\right.
\end{equation}

We call a clockwise (resp. counterclockwise) bubble  fake
if $m+\l-1<0$ and (resp. if $m-\l-1<0$).   These correspond to positive degree bubbles that are labeled by a negative number of dots. These are to be interpreted as formal symbols recursively defined by
 the odd infinite Grassmannian relations
\begin{align} \label{eq:fake-def}
\vcenter{\hbox{$\posbubd{}{2n + \ast \quad \;}$}}
&:=  - \sum\limits_{l=1}^{ n }
\raisebox{2.5mm}{\hbox{$\posbubdfffsl{}{ 2(n-\ell) + \ast \qquad }$}}
\raisebox{-2.5mm}{\hbox{$\negbubdfffsl{}{\; 2l+ \ast }$}}
\qquad \text{for $0 \leq 2n <-\l$,}
\\ \nn
\vcenter{\hbox{$\negbubd{}{\;\;\;2n + \ast }$}}
&:=  - \sum\limits_{l=1}^{ n }
\raisebox{2.5mm}{\hbox{$\posbubdfffsl{}{ 2l+ \ast \;\;}$}}
\raisebox{-2.5mm}{\hbox{$\negbubdfffsl{}{\qquad\;  2(n-\ell) + \ast }$}}
\qquad \text{for $0 \leq 2n <\l$,}
\\ \nn
\vcenter{\hbox{$\posbubd{}{2n+1 + \ast \quad \;}$}}
&:=
\stackrel{\vcenter{\hbox{$\posbubdfffsl{}{ 2n+ \ast \quad }
                    $}}}{\xy (-10,0)*{\;}; (0,0)*{\bigotimes}; \endxy}
\qquad \text{for $0 \leq 2n <-\l$,}
\\ \nn
\vcenter{\hbox{$\negbubd{}{\qquad   2n+1 + \ast }$}}
&:= \quad
\stackrel{\vcenter{\hbox{$\negbubdfffsl{}{ 2n+ \ast  }
                    $}}}{\xy (-4,0)*{\;}; (0,0)*{\bigotimes}; \endxy}
\qquad \text{for $0 \leq 2n +1 <\l$,}
\end{align}

\subsection{$(Q,\Pi)$-envelopes and graded $(Q,\Pi)$-2-supercategories}

\begin{definition}[\cite{BE2} Definition 1.6]
Given a graded 2-supercategory $\mf{U}$, its $(Q,\Pi)$-envelope $\mf{U}_{q,\pi}$ is the graded 2-supercategory with the same objects as $\mf{U}$, 1-morphisms defined from
\[
\Hom_{\mf{U}_{q,\pi}}(\l,\,u) := \left\{
 Q^m\Pi^a F \mid \text{for all $F \in \Hom_{\mf{U}}(\l,\mu)$ with $m\in \Z$ and $a \in \Z_2$}
\right\}
\]
with composition law $(Q^n\Pi^bG)(Q^m\Pi^aF):= Q^{m+n}\Pi^{a+b}(GF)$.  The 2-morphisms are defined by
\[
\Hom_{\mf{U}_{q,\pi}}(Q^m\Pi^aF, Q^n\Pi^bG) :=
\left\{
x_{m,a}^{n,b} \mid \text{for all $x \in \Hom_{\mf{U}}$(F,G)}
\right\}
\]
viewed as a superspace with addition given by $x_{m,a}^{n,b} + y_{m,a}^{n,b} := (x+y)^{n,b}_{m,a}$ and scalar multiplication given by $c(x_{m,a}^{n,b}) := (cx)_{m,a}^{n,b}$. The degrees are given by $\deg( x_{m,a}^{n,b})=\deg(x)+n-m$, $|x_{m,a}^{n,b}|=|x|+a+b$.
The horizontal composition  is given by
\begin{equation} \label{eq:hor-rule}
  y_{m,c}^{n,d}\cdot x_{k,a}^{l,b}
:= (-1)^{c|x|+b|y|+ac+bc} (y \cdot x)_{k+m,a+c}^{l+n,b+d},
\end{equation}
and the vertical composition by
\begin{equation} \label{eq:vert-rule}
 y_{m,b}^{n,b} \circ x_{\ell,a}^{m,b} = (y\circ x)_{\ell,a}^{n,c}.
\end{equation}
\end{definition}

\begin{definition}
We will denote the graded $(Q,\Pi)$-2-supercategory obtained from $(Q,\Pi)$-envelope of $\mf{U}(\mf{sl}_2)$ by $\mf{U}_{q,\pi}=\mf{U}_{q,\pi}(\mf{sl}_2)$.
\end{definition}

\subsection{Idempotent form of categorified quantum group}

 Let $\und{\mf{U}}_{q,\pi}$
denote the underlying 2-category of  $\mf{U}_{q,\pi}$ with the same objects and 1-morphisms but only its even 2-morphisms of degree zero.  Likewise, let $\dot{\und{\mf{U}}}_{q,\pi}$ denote the idempotent completion or Karoubi envelope of $\und{\mf{U}}_{q,\pi}$.  The $(Q,\Pi)$-2-category structure of $\und{\mf{U}}_{q,\pi}$ and $\dot{\und{\mf{U}}}_{q,\pi}$ ensure the existence of distinguished objects $q=(q_{\lambda})$ and $\pi = (\pi_{\lambda})$  in their Drinfeld centers.

We can then define divided powers 1-morphisms $\cal{E}^{(a)}\onel$ and $\cal{F}^{(a)}\onel$ in $\dot{\und{\mf{U}}}_{q,\pi}$.
Recall the idempotent $e_a$ from \eqref{eq:idempotent} in $\End_{\mf{U}}(\onel \cal{E}^a)$ coming from the action of $\ONH_a$.
Define the divided power $\cal{E}^{(a)}\onel$ as
\begin{equation*}
 \cal{E}^{(a)}\onel:=
 \left(Q^{-\frac{a(a-1)}{2}}\cal{E}^a\onel , e_a \right).
\end{equation*}

Let $(e_a')^{{\rm sop}} $ be the element in  $\ONH_a^{{\rm sop}}$ determined from \eqref{eq:idempotent'}. In particular
\[
(e_a')^{{\rm sop}}  :=
(-1)^{\binom{a}{3}+\binom{a}{2}}((x_1^{{\rm sop}})^0 (x_2^{{\rm sop}})^1 \dots (x_{a-1}^{{\rm sop}})^{a-2}(x_a^{{\rm sop}})^{a-1})((\tau_{a-1} ^{{\rm sop}}\dots \tau_1^{{\rm sop}}) (\tau_{a-1} ^{{\rm sop}}\dots \tau_2^{{\rm sop}}) \dots (\tau_{a-1}^{{\rm sop}} \tau_{a-2}^{{\rm sop}})\tau_{a-1}^{{\rm sop}}
\]
Then the divided power $\cal{F}^{(a)}\onel$ is given by
\begin{equation*}
 \cal{F}^{(a)}_i\onel :=
  \left(Q^{\frac{a(a-1)}{2}}\cal{F}^a\onel , (e_a')^{\rm sop} \right).
\end{equation*}
We have that $Q^{-a(\l+a)}\cal{F}^{(a)}\onel$ is right adjoint to $\cal{E}^{(a)}\onel$ in $\dot{\und{\mf{U}}}_{q,\pi}$.
Then by \cite{EKL} there are isomorphisms (recalling Notation~\ref{notation:direct})
 \begin{align*}
  \mathcal{E}^{a}\1_\lambda \cong  [a]_{q,\pi}^!  \mathcal{E}^{(a)}\1_\lambda,
  \qquad \mathcal{F}^{a}\1_\lambda \cong  [a]_{q,\pi}^!  \mathcal{F}^{(a)}\1_\lambda.
\end{align*}

Furthermore, by \cite[Propositions 3.14  and 3.18]{Lau-odd}  there are isomorphisms for $a,b \in \N$
\begin{align}
\label{eq_EaEb}
 \cal{E}^{(a)} \cal{E}^{(b)}\onel  \cong  \qbins{a+b}{a}_{q,\pi}  \cal{E}^{(a+b)}\onel,
 \qquad
 \cal{F}^{(a)}\cal{F}^{(b)}\onel  =  \qbins{a+b}{a}_{q,\pi} \cal{F}^{(a+b)}\onel.
\end{align}
 \begin{align}
\label{eq_FaEb} \cal{F}^{(a)}\cal{E}^{(b)}\onel&=
\bigoplus_{j=0}^{\min(a,b)}  \qbins{a-b-n}{j}_{q,\pi}  \cal{E}^{(b-j)}\cal{F}^{(a-j)}\onel & \text{if $\l < -2a+2$,}\\
\label{eq_EaFb} \cal{E}^{(a)}\cal{F}^{(b)}\onel&=
\bigoplus_{j=0}^{\min(a,b)} \qbins{a-b+n}{j}_{q,\pi} \cal{F}^{(b-j)}\cal{E}^{(a-j)}\onel & \text{if $\l > 2b-2$}.
\end{align}

%

\subsection{2-superfunctors and 2-representations} \label{sec:2-rep}
Recall the notion of a 2-superfunctor  from \cite[Definition 2.2 (ii)]{BE1}
and
  graded $(Q,\Pi)$-2-superfunctor \cite[Definition 5.2 (ii)]{BE2}, see also the equivalent notation of superbifunctor from \cite[Definition 7.2]{KKO2}.
The key point, is that a 2-superfunctor $\Phi \maps \cal{A} \to \cal{B}$ gives superfunctors $\Phi_{X,Y}\maps \Hom_{\cal{A}}(X,Y) \to \Hom_{\cal{B}}(\Phi X, \Phi Y)$ on each Hom supercategory (compatible with the coherence data when $\cal{A}$ and $\cal{B}$ are not strict 2-supercategories).  In addition, if $\Phi$ is a $(Q,\Pi)$-2-superfunctor, then $\Phi_{X,Y}$ is a $(Q,\Pi)$-superfunctor.

For any graded 2-supercategory $\cal{A}$ there is a canonical strict 2-superfunctor $\mathbb{J} \maps \mathcal{K} \to \mathcal{K}_{q,\pi}$ mapping $\cal{A}$ to its $(Q,\Pi)$-envelope $\cal{A}_{q,\pi}$.   If the graded 2-supercategory $\cal{A}$ is already $(Q,\Pi)$ complete, so that there are 1-morphisms $\pi_{\l},q_{\l} \maps \l \to \l$, even 2-isomorphisms $q_{\l} \cong \1_{\l}$, and odd 2-isomorphisms $\pi_{\l} \cong \1_{\l}$ for all objects $\l$ of $\cal{A}$, then $\mathbb{J}$ is a 2-equivalence~\cite[Lemma 4.6]{BE1}.

For a $(Q,\Pi)$-2-supercategory $\cal{B}$, let $\nu \cal{B}$ denote the underlying 2-supercategory.
If $\cal{A}$ is a graded 2-supercategory and $\cal{B}$ is a $(Q,\Pi)$-2-supercategory, then given any 2-superfunctor $\Phi \maps \cal{A} \to \cal{B}$, there is a canonical graded 2-superfunctor $\tilde{\Phi} \maps \cal{A}_{q,\pi} \to \cal{B}$ such that $\Phi=\tilde{\Phi}\mathbb{J}$.  By \cite[Theorem 4.9]{BE1}, there is a functorial equivalence between the 2-category of graded 2-superfunctors $\Hom(\cal{A},\nu \cal{B})$ and the category $\Hom(\cal{A}_{q,\pi},\cal{B})$ of $(Q,\Pi)$-2-superfunctors that sends $\Phi$ to $\tilde{\Phi}$.  In particular, $\tilde{\Phi}$ has the structure of a $(Q,\Pi)$-2-superfunctor.

\begin{definition} \label{def:2rep}
Let $\mf{U}$ be a graded 2-supercategory.
A 2-representation of $\mf{U}$ is a 2-superfunctor $\Phi \maps \mf{U} \to \cal{K}$ for some $(Q,\Pi)$-complete graded 2-supercategory $\cal{K}$.   By the above remarks, such a
 2-representation extends to a canonical $(Q,\Pi)$-2-superfunctor $\tilde{\Phi} \maps \mf{U}_{q,\pi} \to \cal{K}$.
 If in addition, $\cal{K}$ is idempotent complete, then $\Phi$ extends uniquely to a  graded $(Q,\Pi)$-2-superfunctor $\dot{\Phi} \maps \dot{\mf{U}}_{q,\pi} \to \cal{K}$ by the universal property of idempotent completions.

 A morphism of 2-representations $\Phi, \Phi' \maps \mf{U} \to \cal{K}$ is a 2-supernatural transformation $\alpha \maps \Phi \To \Phi'$ in the sense of \cite[Definition 2.2]{BE1}.   Such morphisms extend to $(Q,\Pi)$-2-natural transformations $\tilde{\alpha} \maps \tilde{\Phi} \To \tilde{\Phi'} \maps \mf{U}_{q,\pi} \to \cal{K}$, see \cite[Lemma 6.11]{BE1}
\end{definition}

Breaking down the above definition, given a 2-representation $\Phi \maps \mf{U} \to \cal{K}$, where $\cal{K}$ the 2-supercategory of $(Q,\Pi)$-graded supercategories, then the $(Q,\Pi)$-2-superfunctor $\tilde{\Phi}$ sends each object $\l$ in $\mf{U}_{q,\pi}$ to a $(Q,\Pi)$-supercategory $\Phi(\l)$. Each morphism morphisms of $\mf{U}_{q,\pi}$ maps to $(Q,\Pi)$-superfunctors, and 2-morphisms are sent by $\Phi$ to supernatural transformations.  We will primarily be interested in 2-representations of this form throughout the remainder.

\subsection{2-representations on deformed cyclotomic quotients} \label{sec:2reps-def}

In this section we extend the 2-representations on cyclotomic quotients defined by Kang, Kashiwara, and Oh~\cite{KKO,KKO2}.  The primary difference in our presentation is that the deformed quotients they consider only contain even elements as coefficients, while ours contain odd elements.  More precisely, the $c_i$ are zero in \eqref{eq:an} when $i$ is odd.  Some of the arguments in \cite{KKO} extend in a straight-forward manner, while others are complicated by the additional odd elements.  In this section we focus mainly on those arguments that do not immediately extend.

We start by defining the structure of 2-representation on deformed cyclotomic quotients. The injective morphism
\[
	\left \{ \begin{array}{rcl}
		\ONH_k & \rightarrow & \ONH_{k+1}, \\
		x_a & \mapsto & x_a, \\
		\tau_a & \mapsto & \tau_a, \\
	\end{array} \right.
\]
descends to a morphism $\ONH_k^n \rightarrow \ONH_{k+1}^n$. We denote by $E_n$ (resp. $F_n$) the induction (resp. restriction) superfunctor along these morphisms. We will think of $E_n$ as the $(\ONH_{k+1}^n,\ONH_k^n)$-superbimodule $\ONH_{k+1}^n$. There is an endomorphism $x$ of $E_n$ of degree 2 and parity $\bar 1$ given by right multiplication by $x_{k+1}$. More precisely, $x$ is the $(\ONH_{k+1}^n,\ONH_k^n)$-superbimodule endomorphism defined by
\[
	\funct{\ONH_{k+1}^n}{\ONH_{k+1}^n,}{h}{\varphi(h)x_{k+1},}
\]
with $\varphi$ as in \eqref{eq:varphi-inv}.
Similarly, there is an endomorphism $\tau$ of $E_n^2$ of degree $-2$ and parity $\bar 1$ given by right multiplication by $\tau_{k+1}$. More precisely, $\tau$ is the $(\ONH_{k+2}^n,\ONH_k^n)$-superbimodule endomorphism defined by
\[
	\funct{\ONH_{k+2}^n}{\ONH_{k+2}^n,}{h}{-\varphi(h)\tau_{k+1}.}
\]
Consider the $(Q,\Pi)$-supercategory  of graded supermodules
\[
\cal{L}(n) = \bigoplus_{k=0}^n \sProj(\ONH_k^n) 
\]

\begin{theorem}\label{th:2rep}
	The functors $(E_n,F_n)$ restrict to an adjoint pair of superendofunctors of $\cal{L}(n)$. Together with the data of the natural transformations $x\maps E_n \to \Pi Q^{-2}E_n$ and $\tau\maps E_nE_n \to \Pi Q^2 E_nE_n$,   this defines a 2-representation of $\oUcat_{q,\pi}$ on $\cal{L}(n)$.
\end{theorem}

In the case where $d$ is specialized to 0, this theorem is the type $A_1$ case of \cite[Theorems 8.9 and 9.6]{KKO}. In the rest of this subsection, we explain how to generalize the proof to our case.

The first main argument of the proof is to realize $E_n$ as a cokernel. Consider the following $(\ONH_{k+1},\ONH_k^n)$-superbimodules:
\begin{align*}
	E &= \ONH_{k+1}\otimes_{\ONH_k}\ONH_k^n, \\
	\bar{E} &= \ONH_{k+1,\xi}\otimes_{\ONH_k} \Pi^{n+k}\ONH_k^n,
\end{align*}
where $\ONH_{k+1,\xi} = \ONH_{k+1}$ as a left supermodule, with the right action of $\ONH_k$ given by the morphism
\[
	\xi : \left \{ \begin{array}{rcl}
		\ONH_k & \rightarrow & \ONH_{k+1}, \\
		x_a & \mapsto & x_{a+1}, \\
		\tau_a & \mapsto & \tau_{a+1}. \\
	\end{array} \right.
\]
Then, there is a short exact sequence of $(\ONH_{k+1},\ONH_k^n)$-superbimodules
\begin{eqnarray}\label{eq:ses}
	& & 0 \rightarrow \overline{E} \xrightarrow{P} E \rightarrow E_n \rightarrow 0,
\end{eqnarray}
where the map $E \rightarrow E_n$ is the canonical quotient map, and $P$ is the morphism of degree $2(n-k)$ and parity $\bar{0}$ given by
\[
	\funct{\overline{E}}{E,}{h\otimes\pi^{n+k}y}{ha^n(x_1)\tau_1\ldots\tau_k\otimes y.}
\]
Establishing this exact sequence can be done as in \cite[Section 8]{KKO}. As a consequence, we can deduce that $\ONH_{k+1}^n$ is projective as a right $\ONH_k^n$-supermodule \cite[Theorem 8.7]{KKO}. This proves that $(E_n,F_n)$ send projective modules to projective modules, and thus restrict to superendofunctors of $\cal L(n)$.

A key difference with \cite{KKO}, however, is that the map $P$ is \textit{not} $A_n[t]$-linear, in a sense that we now explain. The superbimodules $\bar{E}, E, E_n$ can be endowed with a structure of right $A_n[t]$-modules given by the following formulas.
\begin{itemize}
	\item On $E_n$, $z\cdot t = zx_{k+1}$.
	\item On $E$, $(h\otimes y)\cdot t = hx_{k+1}\otimes \varphi(y)$.
	\item On $\bar{E}$, $(h\otimes \pi^{n+k}y)\cdot t = (-1)^{k}hx_1\otimes \pi^{n+k}\varphi(y)$.
\end{itemize}
Then, the canonical quotient map $E \twoheadrightarrow E_n$ is $A_n[t]$-linear. But because of relations (\ref{eq:an_relation}), $P$ is not $A_n[t]$-linear. Instead, it is $A_n[t]$-linear \textit{up to an automorphism} of $A_n[t]$. More precisely, we have the relation
\begin{eqnarray}\label{eq:tlin}
	& & P(yt) = P(y)(t-2d).
\end{eqnarray}

Together with Mackey formulas for odd nilHecke algebras (see \cite[Section 5]{KKO}), the short exact sequence (\ref{eq:ses}) gives rise to an exact sequence of $(\ONH_k^n,\ONH_k^n)$-superbimodules
\[
	0 \rightarrow \ker(A) \rightarrow E_n\Pi F_n \xrightarrow{\sigma} F_nE_n \rightarrow \mathrm{coker}(A) \rightarrow 0,
\]
where $A : A_n[t] \otimes_{A_n} \Pi^n\ONH_k^n \rightarrow A_n[t]\otimes_{A_n} \ONH_k^n$ is a certain morphism of $(\ONH_k^n,\ONH_k^n)$-superbimodules induced by $P$ (see \cite[Section 9]{KKO}). In this sequence, the morphism $\sigma$ is the image in $\cal L(n)$ of the 2-morphism 	$\hackcenter{\begin{tikzpicture}[scale=0.4]
		\draw[thick, ->] (0,0) .. controls (0,.5) and (.75,.5) .. (.75,1.0);
		\draw[thick, <-] (.75,0) .. controls (.75,.5) and (0,.5) .. (0,1.0);
		\node at (1.1,.65) {};
	\end{tikzpicture}}$ of the 2-supercategory $\oUcat$.

Contrary to the map $P$, the map $A$ is not $A_n[t]$-linear, even up to an automorphism. However, the default of $A_n[t]$-linearity is controlled. Namely, for all $g(t) \in A_n[t]\otimes_{A_n} Z(\ONH_k^n)$ of degree $\ell$ in $t$, the polynomial $A(yg(t)) - A(y)g((t-2d))$ has degree at most $\ell -1$ in $t$. This can be proved as in \cite[Lemma 9.5]{KKO}, keeping in mind that the $A_n[t]$-linearity of $P$ is replaced by formula (\ref{eq:tlin}) in our case.

Once this is established, we deduce that $A(t^m\otimes \pi^n1)$ is a monic polynomial in $t$ of degree $m+n-2k$, up to a sign. From this, it follows that
\begin{equation*}
	\left \{ \begin{array}{ll}
		\mathrm{ker}(A) = \left(\displaystyle\bigoplus_{i=0}^{n-2k-1}t^i\otimes\ONH_k^n\right) \quad \text{and} \quad \mathrm{coker}(A)= 0 & \text{ if } n-2k \leqslant 0, \\
		\\
		\mathrm{ker}(A) = 0 \quad \text{and} \quad \mathrm{coker}(A) =\left(\displaystyle\bigoplus_{k=0}^{n-2k-1}t^i\otimes \ONH_k^n \right) & \text{ if } n-2k \geqslant0.
	\end{array} \right.
\end{equation*}
The commutation relations \eqref{prop_coveringrelsU1} then follow, completing the proof of Theorem \ref{th:2rep}.

\subsection{Universal quotients and their properties}
 In this section it is convenient to work first with the graded 2-supercategories $\mf{U}$ and its 2-representations and then later pass to its $(Q,\Pi)$-envelope $\mf{U}_{q,\pi}$.

Consider the 2-representation $\Hom_{\oUcat}(-n,-)$ of $\oUcat$. Let $R(n)$ be the 2-subrepresentation generated by $\cal{F}\1_{-n}$. The universal quotient is the 2-representation $\oUcat^n$ of $\oUcat$ defined as
\[
	\oUcat^n = \Hom_{\oUcat}(-n,-) / R(n).
\]
In particular, for all weights $m < -n$, the 1-morphisms $\1_m$ vanish in this quotient. This implies that $\1_m = 0$ for $m>n$.

\begin{proposition}\label{prop:univprop}{\cite[Lemma 5.4]{Rou2}}
	Let $\mathcal V$ be a 2-representation of $\oUcat$, and let $M \in \cal V_{-n}$  be such that $\cal FM=0$. Then there is a canonical morphism of superalgebras $\End_{\oUcat^n}(\1_{-n}) \rightarrow \End_{\mathcal V}(M)$ and a fully-faithful morphism of 2-representations
	\[
	\Phi_M : \oUcat^n\otimes_{\End_{\oUcat^n}(\1_{-n})} \End_{\mathcal V}(M) \rightarrow \mathcal{V}
	\]
	such that $\Phi_M(\1_{-n}) = M$.
	
	Assume furthermore that every object of $\cal V$ is a direct sum of direct summands of objects of the form $\cal E^iM$ for $i \in \N$. Then $\Phi_M$ is an equivalence of 2-representations.
\end{proposition}

\begin{proof}
This proof is analogous to the even case.
\end{proof}

\begin{proposition}\label{prop:univcyclotomic}
	We have an isomorphism of superalgebras $\End_{\oUcat^n}(\1_{-n}) \simeq A_n$ and an equivalence of 2-representations
	\[
		\oUcat^n \simeq \cal L(n)
	\]
	of $\mf{U}$.
	In particular, $\cal{L}(n)$ satisfies the universal property of Proposition \ref{prop:univprop}.  By the discussion in Section~\ref{sec:2-rep}, this equivalence extends to an equivalence $(\oUcat^n)_{q,\pi} \cong \cal{L}(n)$ of 2-representations of $\mf{U}_{q,\pi}$.
\end{proposition}

\begin{proof}
	Put $A_n' = \End_{\oUcat^n}(\1_{-n})$. The object $A_n \in \cal L(n)_{-n}$ satisfies  $F_n(A_n)=0$. Furthermore, every object of $\cal{L}(n)$ is a direct sum of direct summands of objects of the form  $E^i_n(A_n)$. Thus, by the universal property, there is a canonical morphism of superalgebras $A_n' \rightarrow A_n$ and an equivalence of 2-representations
	\[
		\oUcat^n \otimes_{A_n'} A_n \xrightarrow{\sim} \cal L(n)
	\]
	sending $\1_{-n}$ to $A_n$. To prove the proposition, it suffices to show that the canonical morphism $A_n' \rightarrow A_n$ is an isomorphism. We do this by constructing an inverse.

	We know from the proof of the nondegeneracy conjecture~\cite[Theorem 6.2]{DEL} that $\End_{\oUcat}(\1_{-n})$ has a basis given by certain products of bubbles defined in Section~\ref{sec:bubbles}.
	Since $A_n'$ is a quotient of $\End_{\oUcat}(\1_{-n})$, we deduce that (classes of) monomials in bubbles span $A_n'$. However, by definition of $\oUcat^n$, real bubbles vanish in $A_n'$ since they factor through the weight $-n-2$. So $A_n'$ is spanned by monomials in fake bubbles (see \eqref{eq:fake-def}).
	Using relations (\ref{eq:fake-def}), we deduce that there is a surjective morphism
	\[
		\left \{ \begin{array}{rcl}
			A_n & \twoheadrightarrow & A_n', \\[1ex]
			c_i & \mapsto & \hackcenter{ \begin{tikzpicture} [scale=.65]
					\draw[ ]  (0,0) arc (180:360:0.5cm) [thick];
					\draw[<- ](1,0) arc (0:180:0.5cm) [thick];
					\filldraw  [black] (.1,-.25) circle (2.5pt);
					\node at (-.5,-.5) {\tiny $i+\ast$};
					\node at (1.15,.8) {\tiny $-n$};
			\end{tikzpicture}}.
		\end{array} \right.
	\]
where we emphasize that the bubble appearing on the right is a fake bubble, defined by \eqref{eq:fake-def}.
	Let us now prove that the canonical morphism $A_n' \rightarrow A_n$ sends the fake bubble of degree $i$ to $c_i$. The structure of 2-representation of $\cal{L}(n)$ in weights $-n$ and $-n+2$ is given by
\[
\xy
 (-15,0)*+{\sProj(A_n) }="1";
 (25,0)*+{\sProj(A_n[x_1]/a^n(x_1)) }="2";
 {\ar@/^1.3pc/^{E_n = \mathrm{Ind}} "1";"2"};
  {\ar@/^1.3pc/^{F_n = \mathrm{Res}} "2";"1"};
\endxy
\]
	The $\mathfrak{sl}_2$-relations of $\oUcat$ give rise to an isomorphism of $(A_n,A_n)$-superbimodules in $\cal{L}(n)$:
\begin{align}
 \bigoplus_{i=0}^{n-1}
\hackcenter{\begin{tikzpicture}[scale=0.8]
  \draw[thick, ->-=0.15, ->] (-0.37,.5) .. controls ++(-.1,-.8) and ++(.1,-.8) .. (0.35,.5)
      node[pos=0.85, shape=coordinate](Y){}
      node[pos=0.55, shape=coordinate](M){}
      node[pos=0.44, shape=coordinate](X){};
   \node at (Y){\tikz \draw[fill=black] circle (0.4ex);};
\node at (.55,.15) {$\scs i$};
\end{tikzpicture} } \;\maps\;\;
\bigoplus_{i=0}^{n-1}\Pi^i A_n &\to A_n[x_1]/a^n(x_1)  \\
 \nn
(a_0,\ldots,a_{n-1}) &\mapsto \sum_{i=0}^{n-1}x_1^ia_i.
\end{align}
Let $\left[ \alpha_0 \ldots \alpha_{n-1}  \right]$ denote this direct sum of maps and define
	\[
		\left[ \begin{array}{c} \beta_0 \\ \vdots \\ \beta_{n-1} \end{array}  \right] = \left[ \alpha_0 \ldots \alpha_{n-1}  \right]^{-1}
	\]
where $\eta = \alpha_0$ is the unit of the adjunction, see for example \cite[Corollary 5.3]{BE2} for an explicit form of the inverse that utilizes fake bubbles.
	Then by \cite[Definition 2.3]{BE2}, for all $i \in \lbrace 1,\ldots,n\rbrace$ we  can express the fake bubble utilizing these maps as
	\[
		\hackcenter{ \begin{tikzpicture} [scale=.65]
				\draw[ ]  (0,0) arc (180:360:0.5cm) [thick];
				\draw[<- ](1,0) arc (0:180:0.5cm) [thick];
				\filldraw  [black] (.1,-.25) circle (2.5pt);
				\node at (-.5,-.5) {\tiny $i+\ast$};
				\node at (1.15,.8) {\tiny $-n$};
		\end{tikzpicture}} = (-1)^{i+1} \beta_{n-i}\circ x^n \circ \eta.
	\]
	Now, it suffices to observe that
	\begin{align*}
		\left(\beta_{n-i}\circ x^n \circ \eta \right)(1) & = \beta_{n-i}(x_1^n)
		 = \beta_{n-i}\left( \sum_{\ell=0}^{n-1}(-1)^{n-\ell-1}x_1^{\ell}c_{n-\ell}\right)
		 = (-1)^{i-1}c_i.
	\end{align*}
	Hence, the image of the fake bubble of degree $i$ in $A_n$ is $c_i$. It follows that the canonical morphism $A_n' \rightarrow A_n$ is an isomorphism.
\end{proof}

\subsection{Faithfulness of universal quotients}

A 2-representation $\cal V$ of $\oUcat_{q,\pi}$ is said to be \textit{integrable} if the actions of $\cal E, \cal F$ on $\cal V$ are locally nilpotent. A complex $C$ of 1-morphisms of $\oUcat$ is said to be \textit{integrably bounded} if for every integrable representation $\cal V$ and every object $M$ of $\cal V$, the complex $CM$ is bounded.

Define a $(Q,\Pi)$-supercategory
\[
	\cal{L}(n)\mathrm{-sbim} = \bigoplus_{0\leqslant k,\ell \leqslant n} (\ONH_{k}^n,\ONH_{\ell}^n)\mathrm{-sbim}.
\]
The  2-representation from Theorem~\ref{th:2rep}  induces $(Q,\Pi)$-superfunctors $\Phi_{\lambda,n} : \oUcat_{q,\pi}\1_{\lambda} \rightarrow \cal{L}(n)\mathrm{-sbim}$ for all $\lambda \in \mathbb{Z}$. To simplify the notation, we drop the dependence on $\lambda$ and simply denote this superfunctor by
\begin{equation} \label{eq:Phin}
  \Phi_n \maps \oUcat_{q,\pi} \to \cal{L}(n)\mathrm{-sbim}.
\end{equation}

\begin{theorem}\label{th:faithful}
	\begin{enumerate}[(1)]
		\item Let $C$ be an integrably bounded complex of 1-morphisms of $\oUcat$. Assume that for all $n \in \N$, $\Phi_n(C)$ is acyclic. Then, for every integrable 2-representation $\cal V$, $C$ acts by 0 on $K^b(\cal V)$.
		\item Let $f$ be a 2-morphism between integrably bounded complexes of 1-morphisms of $\oUcat$. Assume that for all $n \in \N$, $\Phi_n(f)$ is a quasi-isomorphism. Then, for every integrable 2-representation $\cal V$, $f$ is an isomorphism on $K^b(\cal V)$.
	\end{enumerate}
\end{theorem}

Note that the second point follows from applying the first one to $C=\mathrm{Cone}(f)$. The proof of the first point is decomposed in the following three lemmas. Analogous results were originally proved by Chuang and Rouquier \cite{CR}, \cite{Rou2} in the even case.

\begin{lemma}\label{lem:acyclic}
	Let $C$ be an integrably bounded complex of 1-morphisms of $\oUcat$. Assume that for all $n \in \mathbb{N}$, $\Phi_n(C)$ is acyclic. Then, for every abelian integrable 2-representation $\cal V$ of $\oUcat$ and every object $M$ of $\cal V$, $CM$ is acyclic.
\end{lemma}

\begin{proof}
	We first show the result for $M=\cal E^iN$ where $N$ is a lowest weight object of $\cal V$ and $i \in \N$. Let $n \in \N$  be such that $-n$ is the weight of $N$. By Propositions \ref{prop:univprop} and \ref{prop:univcyclotomic}, we have a morphism of 2-representations $\Phi_N : \cal L(n) \rightarrow \mathcal V$ such that $\Phi_N(A_n)=N$. By assumption, $C\cal E^iA_n$ is acyclic. Since it is a bounded complex of projective supermodules, it is null-homotopic. Furthermore, the complexes $\Phi_N(C\cal E^i A_n)$ and $C\cal E^iN$ are isomorphic since $\Phi_N$ is a morphism of 2-representations. It follows that $C\cal E^iN$ is null-homotopic as well, and therefore acyclic.
	
	We now show the result for an arbitrary object $M$ of $\cal V$. Let $X=C^{\vee}CM$, where $C^{\vee}$ denotes the right dual of $C$. Since $\End_{D^b(\mathcal V)}(CM) \simeq \Hom_{D^b(\mathcal V)}(M,X)$, it suffices to show that $X$ is acyclic. Assume it is not. Denote by $j$ the smallest integer such that $H^j(X) \neq 0$, and let $k$ be the maximal integer such that $\cal F^kH^j(X)\neq 0$. Such integers $j,k$ exist by the boundedness of $X$ and the integrability of $\cal V$ respectively. Put $N=\cal F^kH^j(X)\simeq H^j(\cal F^kX)$, a highest weight object of $\mathcal V$. Then we have isomorphisms
	\[
		\Hom_{D^b(\cal V)}\left(\cal E^kN,X[-j]\right) \simeq \Hom_{D^b(\cal V)}\left(N,\cal F^kX[-j]\right) \neq 0.
	\]
	However, by the first step we know that $C\cal E^kN$ is acyclic, so
	\[
	\Hom_{D^b(\mathcal V)}(\cal E^kN,X[-j]) \simeq \Hom_{D^b(\mathcal V)}(C\cal E^kN,CM[-j])=0,
	\]
	which is a contradiction. Thus $CM$ is acyclic.
\end{proof}

\begin{lemma}
	Let $C$ be an integrably bounded complex of 1-morphisms of $\oUcat$. Assume that for all $n \in \mathbb{N}$, $\Phi_n(C)$ is acyclic. Then for every integrable 2-representation $\cal V$ of $\oUcat$ and every object $M$ of $\cal V$, $CM$ is null-homotopic.
\end{lemma}

\begin{proof}
	We recall the constructions of \cite[5.5.3]{CR}. Consider the object  
	$N =\bigoplus_{i,j\geq 0} \mathcal{F}^i\mathcal{E}^j M$. The full sub supercategory $\cal W$ of $\cal V$ generated by $N$ and closed under direct sums, direct summands, parity and grading shifts, is stable under the actions of $\cal E$ and $\cal F$ by \eqref{eq_EaEb}--\eqref{eq_EaFb}, hence inherits the structure of an integrable 2-representation of $\oUcat$. Let $A = \End_{\cal V}(N)$, a graded superalgebra.  We have an equivalence
	\[
		R = N \otimes_A -: \sProj(A) \xrightarrow{\sim} \cal W
	\]
	such that $M$ is a direct summand of $R(A)$. By transfer of structure, we obtain an integrable 2-representation of $\oUcat$ on $\sProj(A)$ making $R$ into a morphism of 2-representations, see Definition~\ref{def:2rep}. The actions of $\cal E$ and $\cal F$ on $\sProj(A)$ are given by $(A,A)$-superbimodules, thus they also induce a structure of 2-representation of $\oUcat$ on $\sMod(A)$. Applying Lemma \ref{lem:acyclic} to the abelian 2-representation $\mathrm{sMod}(A)$ and the object $A$, we obtain that $CA$ is acyclic. Since it is also a bounded complex of projective $A$-supermodules, we conclude that $CA$ is null-homotopic. Because $R$ is a morphism of 2-representations, there is an isomorphism of complexes $R(CA) \simeq CR(A)$. Thus $CR(A)$ is null-homotopic. Since $M$ is a direct summand of $R(A)$, it follows that $CM$ is null-homotopic as well.
\end{proof}

\begin{lemma}
	Let $C$ be an integrably bounded complex of 1-morphisms of $\oUcat$. Let $\cal V$ be an integrable 2-representation of $\oUcat$ such that for every object $M$ of $\cal V$, $CM$ is null-homotopic. Then $C$ acts by 0 on $K^b(\cal V)$.
\end{lemma}

\begin{proof}
	Observe that $C$ induces a triangulated endofunctor of $K^b(\cal V)$ vanishing on $\cal V$, seen as the full subcategory of $K^b(\cal V)$ consisting of complexes concentrated in cohomological degree 0. Since $\cal V$ generates $K^b(\cal V)$ under cohomological shifts and distinguished triangles, it follows that $C$ vanishes on $K^b(\cal V)$.
\end{proof}

%
\section{Braiding complexes}\label{sec:braid-compl}
%

\subsection{Definition of super equivalences} \label{sec:def-superequiv}

\begin{definition} \label{def:Theta}
  Given $\lambda \in \Z$, we define the \emph{odd Chuang-Rouquier complex}   $\Theta  \onel $  of objects of $\Hom_{\overset{.}{\mf{U}}}(\lambda, -\lambda)$ as follows.
  \begin{itemize}
    \item The $r$'th component of $\Theta  \onel $ is $Q^r\cal F^{(\l+r)}\cal E^{(r)}$.
    \item The differential $d$ is given by the composition of $1_{\cal{F}^{(\l+r)}} \eta 1_{\cal{E}^{(r)}}: \cal{F}^{(\l+r)}\cal{E}^{(r)}\to \cal{F}^{(\l+r)}\cal F\cal E\cal{E}^{(r)} $ with the projection on $\cal{F}^{(\l+r+1)}\cal{E}^{(r+1)}$ given by $\left(e'_{\l+r+1}\right)^{\rm sop}$ and $e_{r+1}$. The fact that the differential squares to zero follows from the fact that $e_2 e_2'=0$.
  \end{itemize}
\end{definition}

 We now consider the images of the odd Chuang-Rouquier complex under the 2-representation $\Phi_n$ from \eqref{eq:Phin}.
To describe the complexes $\Phi_n(\Theta \onel)$ explicitly, we will need additional notation. For $\ell\leqslant m \leqslant k$,  let $e_{[\ell,m]}$ be the idempotent of $\ONH_k$ defined by
\[
	e_{[\ell,m]} = (-1)^{\binom{m-\ell+1}{3}}(\tau_{\ell}\ldots\tau_{m-1})(\tau_{\ell}\ldots\tau_{m-2})\ldots\tau_{\ell}x_{m}^0x_{m-1}^{1}\ldots x_{\ell}^{m-\ell}.
\]
We put $x_{[\ell,m]} = (-1)^{\binom{m-\ell+1}{3}}x_{m}^0x_{m-1}^{1}\ldots x_{\ell}^{m-\ell}$,  and $\tau_{\omega_0[\ell,m]}=(\tau_{\ell}\ldots\tau_{m-1})(\tau_{\ell}\ldots\tau_{m-2})\ldots\tau_{\ell}$, so that $e_{[\ell,m]} = \tau_{\omega_0[\ell,m]}x_{[\ell,m]}$.   Note that if $\ell\leqslant \ell' \leqslant m'\leqslant m$, then $e_{[\ell,m]}e_{[\ell',m']}=e_{[\ell',m']}e_{[\ell,m]}=e_{[\ell,m]}$. We also define an idempotent $e'_{[\ell,m]}$ of $\ONH_k$ by
\[
	e'_{[\ell,m]} = (-1)^{\binom{m-\ell+1}{3}+\binom{m-\ell+1}{2}}x_{\ell}^0x_{\ell+1}^{1}\ldots x_{m}^{m-\ell}(\tau_{m-1}\tau_{m-2}\ldots\tau_{\ell})(\tau_{m-1}\ldots\tau_{2})\ldots\tau_{m-1}.
\]
We put $x'_{[\ell, m]} =(-1)^{\binom{m-\ell+1}{3}+ \binom{m-\ell+1}{2}}x_{\ell}^0x_{\ell+1}^{1}\ldots x_{m}^{m-\ell}$. Note that $ (\tau_{m-1}\tau_{m-2}\ldots\tau_{\ell})(\tau_{m-1}\ldots\tau_{2})\ldots\tau_{m-1} = \tau_{\omega_0[\ell,m]}$ by \cite[(3.51)]{EKL}. Thus, $e'_{[\ell,m]} =x'_{[\ell,m]}\tau_{\omega_0[\ell,m]}$. Note that if $\ell\leqslant \ell' \leqslant m'\leqslant m$, then $e'_{[\ell,m]}e'_{[\ell',m']}=e'_{[\ell',m']}e'_{[\ell,m]}=e'_{[\ell,m]}$.

Let $n \in \N$ and let $k\in \{0,\dots , \frac{n}{2}\}$. Then we have  the complex of $(\ONH_{n-k}^n,\ONH_k^n)$-bimodules
\begin{align*}
	&\Phi_n( \Theta 1_{-n+2k})= \\
	& 0 \to  Q^0\ONH_{n-k}^n e_{[k+1,n-k]} \to \dots \to Q^{r-k}e'_{[n-k+1,r]} \ONH_{r}^n e_{[k+1,r]} \to \dots  Q^{n-k}e'_{[n-k+1,n]} \ONH_{n}^n e_{[k+1,n]} \to 0.
\end{align*}
This complex has $k+1$ non-zero terms, the last one being in cohomological degree $n-k$.

Now let $k\in \{\frac{n}{2}, \dots , n\}$.
Then as complexes of graded $(\ONH_{n-k}^n, \ONH_{k}^n)$-bimodules we have 
\begin{align*}
	&\Phi_n( \Theta 1_{-n+2k}) = \\
	& 0 \to Q^0e'_{[n-k+1,k]} \ONH_{k}^n \to \dots \to Q^{r-k}e'_{[n-k+1,r]} \ONH_{r}^n e_{[k+1,r]} \to \dots  Q^{n-k}e'_{[n-k+1,n]} \ONH_{n}^n e_{[k+1,n]^n} \to 0.
\end{align*}
This complex has $n-k+1$ non-zero terms, the last one being in cohomological degree $n-k$.


\subsection{Graded dimensions and idempotents} \label{sec:grdimidemp}

Suppose that $A$ is a graded superalgebra and $e$ is an idempotent of $A$ of degree and parity 0. Put $B=eAe$, a graded superalgebra. Assume that $Ae$ (resp. $eA$) is free a right (resp. left) $B$-module. If $M$ is a left graded $A$-superbimodule, then $eM$ is naturally a graded $B$-superbimodule, and we have an isomorphism of left graded $A$-supermodules
\[
	\begin{array}{rcl} Ae\otimes_B eM & \xrightarrow{\sim} & M, \\ ae\otimes em & \mapsto & aem. \end{array}
\]
Taking graded dimensions, we get $\grdim(Ae)\grdim(eM)=\grdim(M)\grdim(B)$, from which we deduce
\begin{equation}\label{eq:grdimleftidem}
	\grdim(eM) = \dfrac{\grdim(M)\grdim(eAe)}{\grdim(Ae)}.
\end{equation}
Similarly, if $M$ is a right graded $A$-supermodule, then $Me$ is a right graded $B$-supermodule and
\begin{equation}
	\grdim(Me) = \dfrac{\grdim(M)\grdim(eAe)}{\grdim(eA)}.
\end{equation}

\begin{proposition} \label{prop:eMe}
Given $M$   a left graded $\ONH_m$-supermodule, we have
\begin{equation}\label{eq:grdimlefte'}
	\grdim(e'_mM) = \dfrac{q^{m\choose2}\grdim(M)}{[m]!}, \quad \grdim(e_mM) = \dfrac{q^{-{m\choose2}}\grdim(M)}{[m]!}.
\end{equation}

Likewise, if $M$ is a right graded $\ONH_m$-supermodule, multiplying on the right by the idempotent $e_m$ gives
\begin{equation}\label{eq:grdimrighte}
	\grdim(Me_m) = \dfrac{q^{m\choose2}\grdim(M)}{[m]!}, \quad \grdim(Me'_m) = \dfrac{q^{-{m\choose2}}\grdim(M)}{[m]!}.
\end{equation}
\end{proposition}

\begin{proof}
Let $A=\ONH_{m}$ and $e=e'_{m}$ or $e=e_m$, we have isomorphisms of left graded $\ONH_m$-superbimodules
\[
	\begin{array}{rcl} \OPol_m & \xrightarrow{\sim} & Q^{m(m-1)}\ONH_me'_m, \\ f & \mapsto & f\tau_{\omega_{0}[1,m]}, \end{array} \quad \begin{array}{rcl} \OPol_m & \xrightarrow{\sim} & \ONH_me_m, \\ f & \mapsto & fe_m. \end{array}
\]
as well as isomorphisms of graded superalgebras
\[
	\begin{array}{rcl} \OL_m & \xrightarrow{\sim} & e'_m\ONH_me'_m, \\ f & \mapsto & e'_mfe_m', \end{array} \quad \begin{array}{rcl} \OL_m & \xrightarrow{\sim} & e_m\ONH_me_m, \\ f & \mapsto & e_mfe_m. \end{array}
\]
Hence
\[
	\grdim(\ONH_me'_m) = q^{-m(m-1)}\grdim(\OPol_m) = \dfrac{q^{-m(m-1)}}{(1-q^2)^m},
\]
\[
\grdim(\ONH_me_m) = \grdim(\OPol_m) = \dfrac{1}{(1-q^2)^m},
\]
and
\[
	\grdim(e'_m\ONH_me'_m) = \grdim(e_m\ONH_me_m)= \grdim(\OL_m) = \prod_{i=1}^m \dfrac{1}{1-q^{2i}}.
\]
Therefore
\[
	\dfrac{\grdim(e'_m\ONH_me'_m)}{\grdim(\ONH_me'_m)} = q^{m(m-1)} \prod_{i=1}^m \dfrac{1-q^2}{1-q^{2i}}= q^{m(m-1)} \prod_{i=1}^m \dfrac{1}{q^{i-1}[i]}= \dfrac{q^{m(m-1)}}{q^{m\choose2}[m]!} = \dfrac{q^{m\choose2}}{[m]!},
\]
\[
\dfrac{\grdim(e_m\ONH_me_m)}{\grdim(\ONH_me_m)} = \prod_{i=1}^m \dfrac{1-q^2}{1-q^{2i}}= \prod_{i=1}^m \dfrac{1}{q^{i-1}[i]}= \dfrac{1}{q^{m\choose2}[m]!} = \dfrac{q^{-{m\choose2}}}{[m]!}
\]
and the result follows.
%
\end{proof}

\subsection{Bases}\label{subsubsec:bases}
\begin{definition}
	For any disjoint subsets $U_1,\dots, U_k \subset \{1,\dots, n\}$ of consecutive integers, we write $\OL_n^{U_1,\dots, U_k} \subset \OPol_n$ for the set of odd polynomials that are odd symmetric separately in each of the ordered subsets of consecutive variables $\left \{x_j, \, j\in U_i \right \}$ for $1\leqslant i \leqslant k$.
	
	More precisely, this means that an odd polynomial $f$ is an element of $\OL_n^{U_1,\dots, U_k}$ if and only if $\partial_j(f)=0$ for any $j$ such that $\min(U_i) \leqslant j \leqslant \max(U_i)-1$ for some $1\leqslant i \leqslant k$. For example, $\OL_n^{\emptyset}=\OPol_{n}$, and $\OL_{n}^{\left[1,n\right]} = \OL_n$.
\end{definition}

We fix an integer $n$ and study the action of the odd Chuang-Rouquier complex $\Theta$ in the 2-representation $\cal L(n)$. We start by constructing a basis for $e'_{[\ell,m]}\ONH_m^ne_{[k,m]}$ (which is the form of a general term of $\Phi_n(\Theta)$) as an $A_n$-module.

Define the following sets for $1\leq \ell \leq m$:
\begin{align*}
	X_{m} & = \left \{a=(a_{1},\ldots,a_m), \, a_i \leqslant n-i \right \}, \\
	Y_{\ell,m} & = \left \{a\in X_{m}, \, a_{\ell} > \ldots > a_{m} \right \}.
\end{align*}
We do not indicate the dependence on $n$ (which is fixed for this section) in the notation for simplicity.

\begin{lemma}\label{lem:polbasisME}
	The set $\left \{ \partial_{\omega_0[1,m]}(x^a), \, a \in Y_{1,m}\right \}$ is a basis of $\OL_n^{[1,m],[m+1,n]}$ as a free right ${\OL}_{n}$-module.
\end{lemma}
\begin{proof}
	  For every $a\in Y_{1,m}$, define a partition $\alpha^a=(\alpha^a_1,\dots \alpha^a_m)$ by $\alpha^a_i=a_i-m+i$ for all $i$. This defines a one-to-one correspondence $a\mapsto \alpha^a$ between $Y_{1,m}$ and partitions $\alpha=(\alpha_1,\dots \alpha_m)$ with $\alpha_i\leqslant n-m$ for all $i$. By definition, we have $x^a=\pm x_1^{m-1} \dots x_{m-1} x_m^0 x^{\alpha^a}$. Therefore,
	\[\partial_{\omega_0[1,m]}(x^a)= \pm \partial_{\omega_0[1,m]}(x_1^{m-1} \dots x_{m-1} x_m^0 x^{\alpha^a})=\pm (s_{\alpha^a}(x_1,\dots x_m))^{\omega_{0[1,m]}}\]
 where we use the superscript notation for the action of $w_0[1,m]$   defined below \eqref{eq:sym-act}.
Furthermore, $\OL_n^{[1,m],[m+1,n]}$ has basis as a right $\OL_n$-module given by $s_\l(x_1,\dots x_m)$ for partitions $\lambda=(\l_1, \dots \l_m)$ with $\l_1 \leqslant n-m$ by Lemma~\ref{lem:odd-eqcoh-basis}. Since $\omega_{0[1,m]}$ is an automorphism of $\OL^{[1,m],[m+1,n]}$, the result follows.
\end{proof}

\begin{corollary}\label{cor:5.5}
  $\OL_{n}^{[\ell,n]}$ is a free right $\OL_n$-module with basis given by the set
  \begin{equation}
B_1 =\left\{
x_{1}^{\beta_1}\dots x_{\ell-1}^{\beta_{\ell-1}}s_{\l}(x_1,\dots x_{\ell-1})  \mid \lambda=(\lambda_1,\dots, \lambda_{\ell-1}),  \text{ partition with } \lambda_i \leq n-\ell+1, \; \forall i\;, \beta_i\leq \ell-1-i
\right\} .
\end{equation}
Alternatively, it has a basis given by
  \begin{equation}
B_2 =\left\{
x_{1}^{a_1}\dots x_{\ell-1}^{a_{\ell-1}}  \mid \; a_i\leq n-i
\right\} .
\end{equation}
\end{corollary}
\begin{proof}
  Observe that $\OL_{n}^{[\ell,n]}\cong \OPol_{\ell-1}\otimes_{\OL_{\ell-1}} \OL_{n}^{[1,\ell-1][\ell,n]}$ as a right $\OL_n$-module.
  Then we can multiply the basis $\{x_{1}^{\beta_1}\dots x_{\ell-1}^{\beta_{\ell-1}}\mid \beta_i\leq \ell-1-i\}$ of $\OPol_{\ell-1}$ as a free $\OL_{\ell-1}$-module with the basis $\{s_\l(x_1,\dots x_{\ell-1}) \mid \l=(\l_1,\dots,\l_{\ell-1}) \text{ partition with } \l_i\leq n-\ell+1 \}$ of $\OL_{n}^{[1,\ell-1][\ell,n]}$ as a $\OL_n$-module to conclude that $B_1$ is a basis. To prove that $B_2$ is also a basis, observe that every element in $B_1$ can be written as a $\Z$-linear combination of elements in $B_2$.  Therefore, $B_2$ generates $\OL_{n}^{[\ell,n]}$ as a free $\OL_n$-module. One concludes that $B_2$ is a basis by computing graded dimensions.
\end{proof}

\begin{lemma}\label{lem:polbasisME2}
	The set $\left \{ \partial_{\omega_0[\ell,m]}(x^a), \, a \in Y_{\ell,m}\right \}$ is a basis of $\OL_n^{[\ell,m],[m+1,n]}$ as a right ${\OL}_{n}$-module.
\end{lemma}

\begin{proof}
	Consider the basis of $\OL_{n-\ell+1}^{[1,m-\ell+1],[m-\ell+2,n-\ell+1]}$  as an $\OL_{n-\ell+1}$-module from Lemma \ref{lem:polbasisME}. Tensoring on the left by the algebra $\OPol_{\ell-1}$, we obtain that the set
	\[
		\left \{ \partial_{\omega_0[\ell,m]}(x^a), \, a=(a_{\ell},\ldots,a_m), a_i \leqslant n-i, a_{\ell} >\ldots>a_m \right \}
	\]
	is a basis of $\OL_n^{[\ell,m],[m+1,n]}$ as an $\OL_n^{[\ell,n]}$-module.
 By Corollary \ref{cor:5.5} the set $\left \{ x_1^{a_1}\ldots x_{\ell-1}^{a_{\ell-1}}, a_i\leqslant n-i \right \}$ is a basis of $\OL_n^{[\ell,n]}$ as an $\OL_n$-module.
Multiplying by this basis   of $\OL_n^{[\ell,n]}$ as an $\OL_n$-module
and using the identity
	\[
		x_1^{a_1}\ldots x_{\ell-1}^{a_{\ell-1}}\partial_{\omega_0[\ell,m]}(x^a) = \pm \partial_{\omega_0[\ell,m]}(x_1^{a_1}\ldots x_{\ell-1}^{a_{\ell-1}}x^a)
	\]
	completes the proof.
\end{proof}

\begin{corollary}\label{Cor:Cor5.5ME}
	The set $\left \{ \partial_{\omega_0[\ell,m]}(x^a)\otimes 1, \,  a \in Y_{\ell,m}\right \}$ is a basis of $\OL_m^{[\ell,m]}\otimes_{\OL_m}{M_m^n}$ as a free $A_n$-module.
\end{corollary}
\begin{proof}
We can restrict the isomorphism of $(M_m^n,M_m^n)$-bimodules from Proposition~\ref{prop:eqcoh} to obtain an $(\OL_m, A_n)$-bimodule isomorphism, and since $A_n$ is commutative, we can regard this as an
  isomorphism of $\OL_m\otimes A_n$-supermodules
	\begin{align*}
		& M_m^n \xrightarrow{\sim} (\OL_m\otimes \OL_{n-m})\otimes_{\OL_n} A_n \\
		& (f\otimes c) \mapsto (f\otimes 1)\otimes c.
	\end{align*}
	It induces an isomorphism of $A_n$-modules
	\begin{align*}
		& \OL_m^{[\ell,m]}\otimes_{\OL_m} M_m^n \xrightarrow{\sim} (\OL_m^{[\ell,m]}\otimes \OL_{n-m})\otimes_{\OL_n} A_n \\
		& g\otimes (f\otimes c) \mapsto (gf\otimes 1)\otimes c.
	\end{align*}
	From Lemma \ref{lem:polbasisME2}, the set $\left \{ (\partial_{\omega_0[\ell,m]}(x^a)\otimes 1)\otimes 1, \,  a \in Y_{\ell,m}\right \}$ is a basis of $(\OL_m^{[\ell,m]}\otimes \OL_{n-m})\otimes_{\OL_n} A_n$ as an $A_n$-module. Taking the inverse image of this basis by the isomorphism above yields the result.

\end{proof}


\begin{definition}
	For $a\in Y_{[\ell,m]}$ and $\omega \in \mathfrak{S}_m/\mathfrak{S}_{[k,m]}$, define an element $b_m(a,\omega) \in e'_{[\ell,m]}\ONH_m^ne_{[k,m]}$ by
	\[
	b_m(a,\omega) = e'_{[\ell,m]} x^a\tau_{\omega} e_{[k,m]}\otimes 1.
	\]
\end{definition}

\begin{theorem}\label{th:basis_comp}
	The set  $\{b_m(a,\omega) \mid a\in Y_{\ell,m}, \; \omega\in \mathfrak{S}_m/\mathfrak{S}_{[k,m]}\}$ is a basis of $e'_{[\ell,m]}\ONH_{m}^ne_{[k,m]}$ as an $A_n$-supermodule.
\end{theorem}

\begin{remark}
  Observe that when $k=\ell=m$ this theorem gives a basis for the deformed cyclotomic quotient $\ONH_{m}^n$.
\end{remark}

\begin{proof}
	We start by proving that this set is free over $A_n$. Suppose for contradiction that we have an $A_n$-linear relation $\sum_{a,w} e'_{[\ell,m]}x^a\tau_w e_{[k,m]} \otimes c_{a,w}=0$ in $\ONH_{m}^n$ for $c_{a,w} \in A_n$ not all 0. Let $z \in \mathfrak{S}_m/\mathfrak{S}_{[k,m]}$ be such that $c_{a,z}\neq 0$ for some $a$ and of minimal length for this property.
	
	Consider the action of $\ONH_{m}^n$ on $\OPol_m \otimes_{\OL_m} M_{m}^n$ from Theorem~\ref{thm:matrixMnk}.
 Then consider the action of $\sum_{a,w} e'_{[\ell,m]}x^a\tau_w e_{[k,m]} \otimes c_{a,w}$  on $\mf{s}_z(x_1,\dots x_m)\otimes 1 \in \OPol_m \otimes_{\OL_m} M_{m}^n$, where $\mf{s}_z$ denotes the Schubert polynomial from \eqref{eq-schubert}.
	%
	Observe that
	\begin{align*}
		0 & = \left(\sum_{\substack{a,w\\l(w)\geqslant l(z)}} e'_{[\ell,m]}x^a\tau_w e_{[k,m]}    \otimes c_{a,w}\right) \cdot \mf{s}_{z} \\
		& = \pm \sum_{\substack{a,w\\l(w)\geqslant l(z)}} e'_{[\ell,m]}x^a\partial_w\partial_{z^{-1}w_0[1,m]}(x_1^{m-1}  \dots x_m^0)\otimes c_{a,w} \\
		& = \pm\sum_{\substack{a,w\\l(w)\geqslant l(z)}} e'_{[\ell,m]}x^a\delta_{z,w}\otimes c_{a,w} \\
		& = \pm \sum_a x'_{[\ell,m]}\partial_{\omega_0[\ell,m]}(x^a)\otimes c_{a,z}.
	\end{align*}
	Further applying $\partial_{[\ell,m]}$ to this relation gives
	\[
		\sum_a \partial_{\omega_0[\ell,m]}(x^a)\otimes c_{a,z} = 0
	\]
	From Corollary \ref{Cor:Cor5.5ME}, it follows that $c_{a,z}=0$ for all $z$, a contradiction. Therefore, the linear independence is proved.

To conclude, we compute graded dimensions ignoring parity using Theorem~\ref{thm:matrixMnk} and Proposition~\ref{prop:eqcoh}.
On the one hand, we have
\begin{align*}
  \grdim(\ONH_m^n) &= ([m]!)^2 \grdim(M_m^n) \\
   & =q^{m(n-m)}\frac{[m]![n]!}{[n-m]!}\grdim(A_n)
\end{align*}
Therefore, by (\ref{eq:grdimlefte'}) and (\ref{eq:grdimrighte}), we have
\begin{align*}
  \grdim(e'_{\ell,m}\ONH_{m}^ne_{k,m}) &= \frac{q^{\binom{m-\ell+1}{2}+\binom{m-k+1}{2}}}{[m-\ell+1]![m-k+1]!}\grdim(\ONH_m^n) \\
  &= q^{\binom{m-\ell+1}{2}+\binom{m-k+1}{2}+m(n-m)}\frac{[m]![n]!}{[n-m]![m-\ell+1]![m-k+1]!}\grdim(A_n)
\end{align*}

On the other hand, the free $A_n$-module spanned by the basis $\{e'_{[\ell,m]}x^a\tau_w e_{[k,m]} \mid a\in Y_{\ell,m}, w\in S_{k,m}\}$  has graded dimension
\[
\left(\sum_{a\in Y_{\ell,m}} q^{2(a_1+\dots a_{\ell-1})}q^{2(a_{\ell}+\dots a_m)} \right)\left(\sum_{w\in S_{k,m}} q^{-2\ell(w)}\right)\grdim(A_n)
\]

We have
\begin{align*}
\sum_{a\in Y_{\ell,m}} q^{2(a_1+\dots +a_{\ell-1})}q^{2(a_{\ell}+\dots +a_m)} &=\sum_{a_i\leq n-i} q^{2(a_1+\dots+ a_{\ell-1})}\sum_{\substack{a_i\leq n-i \\a_\ell>a_{\ell+1}>\dots+ >a_m}} q^{2(a_{\ell}+\dots+ a_m)} \\
   &= \sum_{a_i\leq n-i} q^{2(a_1+\dots+ a_{\ell-1})}\sum_{\substack{a_i\leq n-i \\a_\ell>a_{\ell+1}>\dots+ >a_m}} q^{2(a_{\ell}+\dots+ a_m)}  \\
   & = q^{\binom{n}{2}-\binom{n-\ell+1}{2}}\frac{[n]!}{[n-\ell+1]!}\sum_{\substack{a_i\leq n-i \\a_\ell>a_{\ell+1}>\dots+ >a_m}} q^{2(a_{\ell}+\dots+ a_m)} \\
   &= q^{\binom{n}{2}-\binom{n-\ell+1}{2}}\frac{[n]!}{[n-\ell+1]!}q^{2\binom{m-\ell+1}{2}+(n-m)(m-\ell+1)}\frac{[n-\ell+1]!}{[n-m]![m-\ell+1]!}
\end{align*}
with the third equality holding by \cite[(5.7)]{Vera}.

By \cite[(5.6)]{Vera}, we have
\begin{align*}
  \sum_{w\in S_{k,m}} q^{-2\ell(w)} &=  \frac{q^{-\binom{m}{2}} [m]!}{q^{-\binom{m-k+1}{2}} [m-k+1]!}
\end{align*}

Hence the graded dimension as an free $A_n$-module is
$$ q^{\binom{n}{2}-\binom{n-\ell+1}{2}-\binom{m}{2}+\binom{m-k+1}{2}+2\binom{m-\ell+1}{2}+(n-m)(m-\ell+1)}\frac{[n]![m]!}{[n-m]![m-\ell+1]![m-k+1]!}$$
The result follows since the power of $q$ is equal to $\binom{m-k+1}{2}+\binom{m-\ell+1}{2}+m(n-m)$.

	\end{proof}

\subsection{Invertibility}

\subsubsection{Differential}

In the general form considered for the terms of $\Phi_n(\Theta)$, the differential of $\Phi_n(\Theta)$ takes the form
\[
	{\sf d}_m = \left \{ \begin{array}{rcl}
		e'_{[\ell,m]} \ONH_m^n e_{[k,m]} & \rightarrow & Qe'_{[\ell,m+1]}\ONH_{m+1}^ne_{[k,m+1]}, \\
		h & \mapsto & e'_{[\ell,m+1]}he_{[k,m+1]}.
	\end{array}\right.
\]

\begin{theorem}\label{th:rickard_exact}
	We have $\ker({\sf d}_{m+1})=\mathrm{im}({\sf d}_m)$.
\end{theorem}

\begin{proof}
	We start by computing ${\sf d}_m(b_m(a,\omega)) = e'_{[\ell,m+1]}x^a\tau_{\omega}e_{[k,m+1]}$ for $a\in Y_{\ell,m}$ and $\omega \in \mf{S}_{m}/\mf{S}_{[k,m]}$.\\	
	If $a_m >0$, then we have ${\sf d}_m(b_m(a,\omega))=b_{m+1}((a,0),\omega)$, where $(a,0)=(a_1,\ldots,a_m,0) \in Y_{\ell,m+1}$ and $w$ is seen as an element of $\mf{S}_{m+1}/\mf{S}_{[k,m+1]}$ fixing $m+1$. \\
	In the case $a_m=0$, let $r$ be the largest integer with $r\leqslant m-\ell$ such that $a_{m}=0<a_{m-1}=1<\ldots<a_{m-r}=r$. The key point of the computation is to observe that we have the following equality in $\ONH_{m+1}$:
	\[
	\tau_{m-r}x_{m-r+1}^{r+1}x_{m-r}^r + x_{m-r}^{r+1}x_{m-r+1}^{r}\tau_{m-r} = x_{m-r}^{r}x_{m-r+1}^r
	\]
	It follows that
	\begin{align*}
		& \tau_{\omega_0[\ell,m+1]}x_{\ell}^{a_{\ell}}\ldots x_{m-r+1}^{a_{m-r+1}}x_{m-r}^{r+1}x_{m-r+1}^r\ldots x_{m}\tau_{m-r} \\
		= \ & (-1)^{\frac{r(r-1)}{2}} \tau_{\omega_0[\ell,m+1]}x_{\ell}^{a_{\ell}}\ldots x_{m-r+1}^{a_{m-r+1}}\left(x_{m-r}^{r+1}x_{m-r+1}^r\tau_{m-r}\right)x_{m-r+2}^{r-1}\ldots x_{m} \\
		= \ & (-1)^{\frac{r(r-1)}{2}} \tau_{\omega_0[\ell,m+1]}x_{\ell}^{a_{\ell}}\ldots x_{m-r+1}^{a_{m-r+1}}\left(x_{m-r}^{r}x_{m-r+1}^r\right)x_{m-r+2}^{r-1}\ldots x_{m} \\
		\ & - (-1)^{\frac{r(r-1)}{2}} \tau_{\omega_0[\ell,m+1]}x_{\ell}^{a_{\ell}}\ldots x_{m-r+1}^{a_{m-r+1}}\left(\tau_{m-r}x_{m-r+1}^{r+1}x_{m-r}^r\right)x_{m-r+2}^{r-1}\ldots x_{m}. \\
	\end{align*}
	This last term is $0$ since $\tau_{\omega_0[\ell,m+1]}\tau_{m-r}=0$. Hence, we have
	\begin{align*}
		& \tau_{\omega_0[\ell,m+1]}x_{\ell}^{a_{\ell}}\ldots x_{m-r+1}^{a_{m-r+1}}x_{m-r}^{r+1}x_{m-r+1}^r\ldots x_{m}\tau_{m-r} \\
		= \ & (-1)^{\frac{r(r-1)}{2}} \tau_{\omega_0[\ell,m+1]}x_{\ell}^{a_{\ell}}\ldots x_{m-r+1}^{a_{m-r+1}}x_{m-r}^{r}x_{m-r+1}^rx_{m-r+2}^{r-1}\ldots x_{m}
	\end{align*}
	We can repeat this argument for $\tau_{m-r+1},\ldots,\tau_{m}$ and we obtain
	\begin{align*}
		& \tau_{\omega_0[\ell,m+1]}x_{\ell}^{a_{\ell}}\ldots x_{m-r+1}^{a_{m-r+1}}x_{m-r}^{r+1}x_{m-r+1}^r\ldots x_{m}\tau_{m-r}\ldots\tau_m \\
		= & (-1)^{\frac{r(r-1)(r+1)}{6}} \tau_{\omega_0[\ell,m+1]} x_{\ell}^{a_{\ell}}\ldots x_{m-r+1}^{a_{m-r+1}}x_{m-r}^{r}x_{m-r+1}^{r-1}\ldots x_{m-1}.
	\end{align*}
	It follows that
	\begin{equation}\label{eq:diffbasis}
		{\sf d}_m(b_m(a,\omega)) = (-1)^{\frac{r(r-1)(r+1)}{6}}e'_{[\ell,m+1]}x_{1}^{a_{1}}\ldots x_{m-r+1}^{a_{m-r+1}}x_{m-r}^{r+1}x_{m-r+1}^r\ldots x_{m}\tau_{m-r}\ldots\tau_m\tau_{\omega}e_{[k,m+1]}.
	\end{equation}
	Observe that if $w(m)\geqslant m-r$, then $s_{m-r}\ldots s_mw$ is not an element of minimal length in its left coset modulo $\mf{S}_{[k,m+1]}$, and therefore $\tau_{m-r}\ldots\tau_m\tau_{\omega}e_{[k,m+1]}=0$. Otherwise if $w(m)<m-r$, then $\tau_{m-r}\ldots\tau_m\tau_{\omega}e_{[k,m+1]} = \tau_{s_{m-r}\ldots s_m\omega}e_{[k,m+1]}$, where $s_{m-r}\ldots s_m\omega$ is a well-defined element of $\mf{S}_{m+1}/\mf{S}_{[k,m+1]}$.\\
	Hence, we have shown that $\d_m(b_m(a,w))$ can be computed as follows.
	\begin{itemize}
		\item 
For $r$ as above and $w(m)\geqslant m-r$, then
		\[
			\d_m(b_m(a,w))=0.
		\]
		\item 
For $r$ as above and $w(m)< m-r$, then
		\[
			\d_m(b_m(a,w))=\pm b_{m+1}((a_1,\ldots,a_{m-r+1},r+1,\ldots,0),s_{m-r}\ldots s_mw).
		\]
	\end{itemize}

	We can now deduce bases of $\ker(\d_m)$, $\im(\d_m)$ from these computations.

For $0\leq r\leq m-\ell$, define \[Y_{\ell,m}^r:=\{(a_1,\dots, a_m)\in Y_{\ell,m} \mid \forall i\leq r, \; a_{m-i}=i \text{ and } a_{m-r-1}>r+1 \text{ if } r<m-\ell \}\]
For $1\leq u \leq m$, define $S_{k,m}^u$ (resp. $S_{k,m}^{\geq u}$ $S_{k,m}^{<u}$ ) as the subset of $\mf{S}_{m}/\mf{S}_{[k,m]}$ consisting of elements $w$ with $w(m)=u$ (resp. $w(m)\geq u$, $w(m)<u$).
Then we have proven that the sets
  \[
  \left\{ b_m(a,w) \mid  (a,w) \in \bigsqcup\limits_{r=0}^{m-\ell} Y_{\ell,m}^r\times S_{k,m}^{\geq m-r} \right\}
  \]
  \[
  \left\{ b_{m+1}(a,w) \mid  (a,w) \in \bigsqcup\limits_{r=0}^{m-\ell} Y_{\ell,m+1}^r\times S_{k,m+1}^{\geq m+1-r} \right\}
  \]
  are bases for $\ker(\d_m)$ and $\im(\d_m)$ respectively as free $A_n$-modules.
This proves the theorem.
\end{proof}

\begin{corollary}
  For all $n\in\N$, $k\in \{0, \dots, n\}$, and $i\neq n-k$, we have:
  \[ H^{i}\left(\Phi_n(\Theta\1_{-n+2k})\right)=0.\]
\end{corollary}

\subsubsection{Top Cohomology}

As a consequence of Theorem \ref{th:rickard_exact}, we deduce the following corollary.

\begin{corollary}
	For all $\lambda \in \mathbb Z$, the cohomology of $\Phi_n(\Theta\1_{\lambda})$ is concentrated in top cohomological degree.
\end{corollary}

We now prove that the top cohomology of $\Phi_n(\Theta\1_{-n+2k})$ is an invertible $(\ONH_{n-k}^n, \ONH_k^n)$-superbimodule. Put
\[
	C_{n,k} = H^{\bullet}\left(\Phi_n(\Theta\1_{-n+2k})\right).
\]

We start by giving a basis of $C_{n,k}$ over $A_n$.
\begin{lemma}\label{lem:cnk_basis}
	As an $A_n$-module, $C_{n,k}$ is free with basis
	$$\{x_1^{a_1}\dots x_{n-k}^{a_{n-k}} e'_{[n-k+1,n]} x_{[n-k+1,n]} \tau_{\sigma_k} \tau_w e_{[k+1,n]}\otimes 1 \mid a_i\leq n-i, \; w\in \mf{S}_k \},$$
	where $\sigma_k$ the longest element of $\mf{S}_n/\mf{S}_k\times \mf{S}_{n-k}$.
\end{lemma}

\begin{proof}
	From Theorem \ref{th:basis_comp}, we know that a basis of $e'_{[n-k+1,n]}\ONH_{n}^ne_{[k+1,n]}$ as an $A_n$-module is given by
	\[
		\left \{ e'_{[n-k+1,n]}x^a\tau_we_{[1,k]}, a \in Y_{n-k+1,n}, w\in  \mf{S}_n/\mf{S}_{[k+1,n]}\right \}.
	\]
	From Equation \ref{eq:diffbasis}, we know that a basis of $\im(\d_{n-1})$ as an $A_n$-module is given by
	\[
	\left \{ e'_{[n-k+1,n]}x^a\tau_we_{[1,k]}, a \in Y_{n-k+1,n}, w\in  \mf{S}_n/\mf{S}_{[k+1,n]}, w(n)\geq n-k+1\right \}.
	\]
	Therefore, a basis for $C_{n,k}$ as an $A_n$-module is given by
	\[
	\left \{ e'_{[n-k+1,n]}x^a\tau_we_{[1,k]}, a \in Y_{n-k+1,n}, w\in  \mf{S}_n/\mf{S}_{[k+1,n]}, w(n)< n-k+1\right \}.
	\]
	Note that if $a\in Y_{n-k+1,n}$, then $a=(a_1,\ldots,a_{n-k},k-1,\ldots,0)$ with $a_i\leqslant n-i$. Furthermore, any $w\in \mf{S}_n/\mf{S}_{[k+1,n]}$ such that $w(n)< n-k+1$ can be written as $\sigma_kw'$ with $w'\in \mf{S}_k$. The result follows.
\end{proof}

\begin{corollary}\label{lem:ecnke_basis}
	The following set forms a basis for $e'_{[1,n-k]} C_{n,k} e_{[1,k]}$ as an $A_n$-module.
$$\{x'_{[1,n-k]} \partial_{w_0[1,n-k]}(x_1^{a_1}\dots x_{n-k}^{a_{n-k}}) x'_{[n-k+1,n]} \tau_{w_0[1,n]} x_{[k+1,n]}x_{[1,k]}\otimes 1 \mid n-1\geqslant a_1>a_2>\ldots >a_{n-k} \geqslant 0 \}.$$
\end{corollary}

\begin{proof}
	We start by proving that this set spans $e'_{[1,n-k]} C_{n,k} e_{[1,k]}$. To do this, we show that if $b$ is an element of the basis from Lemma \ref{lem:cnk_basis}, then $e'_{[1,n-k]}be_{[1,k]}$ is in the span of this set. First, notice that for $w \in \mf{S}_k$, we have
	\[
		\tau_we_{[k+1,n]}e_{[1,k]} = \begin{cases}
			\pm \tau_{\omega_{0}[1,k]}\tau_{\omega_{0}[k+1,n]}x_{[1,k]}x_{[k+1,n]} & \text{ if } w=1,\\
			0 & \text{ if } w\neq 1.
		\end{cases}
	\]
	Furthermore, for all $(a_1,\ldots,a_{n-k})$ we have
	\begin{align*}
		& e'_{[1,n-k]}x_1^{a_1}\ldots x_{n-k}^{a_{n-k}}e'_{[n-k+1,n]}x_{[n-k+1,n]}\tau_{\sigma_k}\tau_{\omega_{0}[1,k]}\tau_{\omega_{0}[k+1,n]}\\
		= \ & \pm e'_{[1,n-k]}x_1^{a_1}\ldots x_{n-k}^{a_{n-k}}e'_{[n-k+1,n]}x_{[n-k+1,n]}\tau_{\omega_{0[1,n]}} \\
		= \ & \pm x'_{[1,n-k]} \partial_{w_0[1,n-k]}(x_1^{a_1}\dots x_{n-k}^{a_{n-k}}) x'_{[n-k+1,n]} \tau_{w_0[1,n]}.
	\end{align*}
	It follows that the given set spans $e'_{[1,n-k]} C_{n,k} e_{[1,k]}$ as an $A_n$-module. To conclude that it is a basis, it suffices to compute graded dimensions.  By Proposition~\ref{prop:eMe} we have that
\[
\grdim (e'_{[1,n-k]} C_{n,k} e_{[1,k]} )
= \frac{q^{\binom{n-k}{2}}}{[n-k]!}\frac{q^{\binom{k}{2}}}{[k]!} \grdim ( C_{n,k} ) .
\]

\end{proof}

Since each chain groups of $\Phi_n(\Theta\1_{\lambda})$ and its top cohomology $C_{n,k}$ are $(\ONH_{n-k}^n,\ONH_k^n)$-bimodules with left and right action given by left and right multiplication respectively, it follows from the discussion in Section~\ref{sec:grdimidemp} that $e'_{[1,n-k]} C_{n,k}e_{[1,k]}$ is a $(M_{n-k}^n, M_k^n)$-bimodule with left action of $f\in M_{n-k}^n$ given by left multiplication by $e'_{[1,n-k]}f$ and right action of $g\in M_{n-k}^n$ given by right multiplication by $ge_{[1,k]}$. We now study this bimodule and prove its invertibility.

Put $b_{k,n}=x'_{[1,n-k]}x'_{[n-k+1,n]} \tau_{w_0[1,n]} x_{[k+1,n]}x_{[1,k]}$. This is an element of the basis of $e'_{[1,n-k]} C_{n,k}e_{[1,k]}$ given in Corollary \ref{lem:ecnke_basis}. It has  degree $2((n-k)(n-k-1)+k(k-1))-n(n-1)$ and parity $\binom{n}{2}$.

\begin{lemma}\label{lem:left_action}
	As a left $M_{n-k}^n$-module, $e'_{[1,n-k]} C_{n,k}e_{[1,k]}$ is free of rank 1 with basis $\left \{ b_{k,n} \right \}$.
\end{lemma}

\begin{proof}
	Consider the morphism of left $M_{n-k}^n$-supermodules
	\[
	\begin{array}{rcl}
		\Omega\maps \Pi^{\binom{n}{2}}M_{n-k}^n & \rightarrow & e'_{[1,n-k]} C_{n,k}e_{[1,k]}, \\
		\pi^{\binom{n}{2}}f\otimes a  & \mapsto & \varphi^{\binom{n}{2}}(f\otimes a)\cdot b_{k,n}.
	\end{array}
	\]
	For a partition $\lambda = (\lambda_1, \dots, \lambda_{n-k})$, we have
	\begin{align*}
		\Omega\left(\pi^{\binom{n}{2}}s_{\lambda}\otimes 1\right) & = \varphi^{\binom{n}{2}}(s_{\lambda}\otimes 1)\cdot b_{k,n} \\
		& = \pm x'_{[1,n-k]}s_\lambda^{w_0[1,n-k]} x'_{[n+1-k,n]}\tau_{\omega_0[1,n]}x_{[1,k]}x_{[k+1,n]}.
	\end{align*}
	It follows from Corollary~\ref{lem:ecnke_basis} that $\Omega$ sends a basis to a basis, and therefore is an isomorphism. This proves the lemma.
\end{proof}

\begin{theorem}\label{th:cnkiso}
	The $(M_{n-k}^n,M_k^n)$-superbimodule $e'_{[1,n-k]} C_{n,k}e_{[1,k]}$ is invertible.
\end{theorem}

\begin{proof}
	Since $\left \{ b_{k,n} \right \}$ is a basis of $e'_{[1,n-k]} C_{n,k}e_{[1,k]}$ as a left $M_{n-k}^n$-module by Lemma \ref{lem:left_action}, there exists a unique morphism of superalgebras $u:M_k^n\rightarrow M_{n-k}^n$ defined by
	\[
		b_{k,n}\cdot m = u(m)\cdot b_{k,n}
	\]
	for all $m\in M_{k}^n$. To prove the invertibility of $e'_{[1,n-k]} C_{n,k}e_{[1,k]}$ as an $(M_{n-k}^n,M_k^n)$-superbimodule, it suffices to show that $u$ is an isomorphism. Since $\grdim(M_k^n)=\grdim(M_{n-k}^n)$ by Proposition \ref{lem:odd-eqcoh-basis}, it is enough to show that $u$ is injective, which we can do by reducing modulo 2 by Corollary \ref{cor:inj_mod2}.
	
	Working modulo 2, we have that $\OL_n\otimes_{\Z}\Z_2$ is a central subalgebra of the commutative algebra $\OPol_{n}\otimes_{\Z}\Z_2$ and the algebra $\ONH_n\otimes_{\Z}\Z_2$. Thus, if $p\in \OL_n$, then as elements of $\ONH_n^n\otimes_{\Z} \Z_2$, we have
	\begin{equation}\label{eq:central}
		(p\otimes 1)\otimes 1=(1\otimes \bar{p})\otimes 1
	\end{equation}
	where $\bar{p}$ is the image of $p$ in $A_n$ via the quotient map (\ref{eq:an_quotient_map}).
This centrality also allows to endow $\OL_n^{[1,k],[k+1,n]}\otimes_{\OL_n}A_n\otimes_{\Z}\Z_2$ with an algebra structure via (\ref{eq:superalg-tensor}) making the map
	\[
		\beta_k^n\otimes_{\Z} \Z_2 : M_k^n \otimes_{\Z} \Z_2 \rightarrow \underline{\OL}_n^{[1,k],[k+1,n]}\otimes_{\Z}\Z_2
	\]
	induced from the isomorphism of Proposition \ref{prop:eqcoh} an isomorphism of algebras,   where for notational simplicity, given a right-module over $\OL_n$, we will write $\underline{M}:= M\otimes_{\OL_n} A_n$.  From (\ref{eq:central}), it follows that the right (resp. left) action of $M_k^n \otimes_{\Z} \Z_2$ (resp. $M_{n-k}^n \otimes_{\Z} \Z_2$) on $\ONH_n^n\otimes_{\Z}\Z_2$ satisfies the property
	\begin{align*}
		& (h\otimes 1\otimes 1)\cdot (\beta_k^n\otimes_{\Z} \Z_2)^{-1}(f\otimes 1\otimes 1) = hf\otimes 1\otimes 1\\
		& \left(\text{resp. } (\beta_{n-k}^n\otimes_{\Z} \Z_2)^{-1}(f\otimes 1\otimes 1)\cdot(h\otimes 1\otimes 1) = fh\otimes 1\otimes 1\right)
	\end{align*}
	for $h\in \ONH_n$ and $f\in \OL_n^{[1,k],[k+1,n]}$ (resp. $f\in \OL_n^{[1,n-k],[n-k+1,n]}$).
	
	We now compute the map $u\otimes_{\Z}\Z_2$. Let $f\in\OL_n^{[1,k],[k+1,n]}$. Then in $e'_{[1,n-k]} C_{n,k}e_{[1,k]}\otimes_{\Z} \Z_2$ we have
	\begin{align*}
		& (b_{k,n}\otimes 1) \cdot (\beta_k^n\otimes_{\Z} \Z_2)^{-1}(f\otimes 1\otimes 1) \\
		= \ & x'_{[1,n-k]}x'_{[n-k+1,n]} \tau_{w_0[1,n]} x_{[k+1,n]}x_{[1,k]}f\otimes 1\otimes 1 \\
		= \ & \omega_0(f)x'_{[1,n-k]}x'_{[n-k+1,n]} \tau_{w_0[1,n]} x_{[k+1,n]}x_{[1,k]}\otimes 1\otimes 1 \ \text{ modulo coboundaries} \\
	 	= \ & (\beta_{n-k}^n\otimes_{\Z}\Z_2)^{-1}(\omega_{0}(f)\otimes 1\otimes 1)\cdot (b_{k,n}\otimes 1).
	\end{align*}

	So we have proved that the diagram of algebras
	\[
		\xymatrix{
		M_k^n\otimes_{\Z}\Z_2 \ar[r]^-{u\otimes_{\Z}\Z_2} \ar[d]^-{\sim}_-{\beta_k^n\otimes_{\Z}\Z_2} & M_{n-k}^n\otimes_{\Z}\Z_2 \ar[d]^-{\beta_{n-k}^n\otimes_{\Z}\Z_2}_-{\sim}\\
		\underline{\OL}_n^{[1,k],[k+1,n]}\otimes_{\Z}\Z_2 \ar[r]_-{\omega_{0}}^-{\sim} & \underline{\OL}_n^{[1,n-k],[n-k+1,n]}\otimes_{\Z}\Z_2 \\
		}
	\]
	commutes. The result follows.
\end{proof}

	The following corollary follows immediately from Theorem \ref{th:cnkiso}.
	\begin{corollary}\label{cor:invert_bimod}
		For all $k\leqslant n$, $H^{\bullet}\left(\Phi_n(\Theta\1_{-n+2k})\right)$ is an invertible $\left(\ONH_{n-k}^n,\ONH_{k}^n\right)$-bimodule.
	\end{corollary}

	\begin{theorem}\label{thm:main}
		Let $\cal V$ be an integrable 2-representation of $\oUcat$. For all $\lambda \in \Z$, $\Theta \1_{\lambda}$ induces an equivalence $K^b(\cal V_{\lambda}) \xrightarrow{\sim} K^b(\cal V_{-\lambda})$.
	\end{theorem}

	\begin{proof}
		Denote by $\Theta^{\vee}\1_{-\lambda}$ the right adjoint of $\Theta \1_{\lambda}$, so that $\Theta^{\vee}\1_{-\lambda}$ is a complex of 1-morphisms $-\lambda \rightarrow \lambda$. Consider the unit of adjunction $\eta : \1_{\lambda} \rightarrow \Theta^{\vee}\Theta \1_{\lambda}$. By Corollary \ref{cor:invert_bimod}, $\Phi_n(\eta)$ is a quasi-isomorphism for all $n \in \N$. By Theorem \ref{th:faithful}, it follows that $\eta$ is an isomorphism on $K^b(\cal V_{\lambda})$. The same argument can be applied to the counit of adjunction $\varepsilon : \Theta \Theta^{\vee}\1_{-\lambda} \rightarrow \1_{-\lambda}$, proving that it is an isomorphism on $K^b(\cal V_{-\lambda})$. The result follows.
	\end{proof}
%
\section{Applications }
%

\subsection{Quiver Hecke superalgebras} \label{sec:QHA}

Let $I$ be an index set partitioned into $I=I_{{\rm even}} \sqcup I_{{\rm odd}}$ and $A=(a_{ij})_{i,j \in I}$ a symmetrizable Cartan matrix satisfying $a_{ij}\in 2\Z$ for all $i\in I_{{\rm odd}}$ and $j\in I$.  Define a parity function $p\maps I \to \{0,1\}$ by $p(i)=1$ if $i \in I_{{\rm odd}}$ and $p(i)=0$ if $i\in I_{{\rm even}}$.
Let $\cal{P}_{i,j}:=\Bbbk\la w,z\ra / \la zw- (-1)^{p(i)p(j)}wz\ra$ be the $\Z\times \Z_2$-graded $\Bbbk$-algebra where $w$ and $z$ where $w$ and $z$ have $\Z\times \Z_2$-degree $((\alpha_i,\alpha_i, p(i) )$ and $((\alpha_j,\alpha_j, p(j) )$, respectively.


To a Cartan datam and a matrix $\cal{Q}_{i,j} \in \cal{P}_{i,j}$ satisfying various conditions, Kang, Kashiwara, and Tsuchioka define a \emph{quiver Hecke superalgebra} $R_n = R_n(\cal{Q})$.  The odd nilHecke algebra studied in Section~\ref{sec:oddnil} is the rank one quiver Hecke superalgebra with a single odd $i\in I$.  The quiver Hecke algebra $R_n$ decomposes into blocks $R_{\beta}$ for $\beta=\sum_u m_i\alpha_i\in Q_+$ with $\sum_im_i=n$.  For dominant integral weight $\Lambda \in \P_+$, the quiver Hecke superalgebra admits a cyclotomic quotient $R^{\Lambda}_n$.  The direct sum $R^{\Lambda}:=\oplus_n R^{\Lambda}_n$ categorifies the representation $V(\Lambda)$ of the Kac-Moody superalgebra associated with the Cartan data~\cite{KKO,KKO2}.

The quiver Hecke superalgebra extends to a 2-supercategory $\mf{U}(\mf{g})$ introduced by Brundan and Ellis~\cite{BE2}.   For a single odd vertex, $i \in I$, this 2-supercategory agrees with the 2-supercategory $\mf{U}(\mf{sl}_2)$ from Definition~\ref{def:oddU}.

\subsection{2-representations}

Let $A$ be a $\Z$-graded superalgebra.
To avoid repetition, let $X\in \{\emptyset, \rm{super} \}$ so that $\xMod(A)$ can be used to denote either the supercategory $\Mod(A)$ of left $A$ modules if $X=\emptyset$, or the supercategory $\sMod(A)$ of left $A$ supermodules with $\Z_2$-degree preserving homomorphisms if $X= \rm{super}$.  Likewise, we denote by $\xRep(A)$ the supercategory of $\Z$-graded $A$-modules (resp. $A$-supermodules) that are finite-dimensional over $\Bbbk_0$ if $X=\emptyset$ (resp. $X= \rm{super}$).  Then $\xProj(A)$ denotes either the supercategory of finitely generated projective $A$-modules if $X=\emptyset$, or the supercategory of finitely generated projective $A$-supermodules if $X=\rm{super}$.

By the results in Section 9 of \cite{KKO} for $X=\emptyset$, and
by Theorem 8.13 \cite{KKO2} and the results in Section 8.3 for $X=\rm{super}$, we have the following result.

\begin{proposition}
For each $i\in I$ and $\beta \in Q^+$, there exist $(Q,\Pi)$-superfunctors
\begin{align}
  E_i^{\Lambda} & \maps \xMod(R^{\Lambda}(\beta+\alpha_i)) \to \xMod(R^{\Lambda}(\beta)) \nn \\
  F_i^{\Lambda} & \maps \xMod(R^{\Lambda}(\beta)) \to \xMod(R^{\Lambda}(\beta+\alpha_i))
\end{align}
that are exact on $\xRep(R^{\Lambda})$ and $\xProj(R^{\Lambda})$.  Furthermore, there exits natural isomorphisms of endofunctors on $\xMod(R^{\Lambda}(\beta))$ given by
\begin{align}
&  E_i^{\Lambda}  F_j^{\Lambda} \to q^{-(\alpha_i \mid \alpha_j)} \Pi^{p(i)p(j)}  F_j^{\Lambda}E_i^{\Lambda}  \quad \text{if $i\neq j$,}
\\
&  \Pi_i q_i^{-2} F_i^{\Lambda}E_i^{\Lambda}  \oplus \bigoplus_{k=0}^{\langle h_i, \Lambda-\beta \rangle -1} \Pi_i^k q_i^{2k}
  \to  E_i^{\Lambda}  F_i^{\Lambda} \quad \text{if $\langle h_i, \Lambda-\beta\rangle \geq 0$ },
   \\
&   \Pi_i q_i^{-2} F_i^{\Lambda}E_i^{\Lambda}  \to E_i^{\Lambda}  F_i^{\Lambda}
    \oplus \bigoplus_{k=0}^{-\langle h_i, \Lambda-\beta \rangle -1} \Pi_i^{k+1} q_i^{-2k-2}
   \to  \quad \text{if $\langle h_i, \Lambda-\beta\rangle < 0$},
\end{align}
Then the main result of Brundan and Ellis~\cite{BE2} implies that there is a 2-representation of $\mathfrak{U}(\mf{g})$   on the supercategories
\[
\xRep(R^{\Lambda}) := \bigoplus_{\beta \in Q^+}  \xRep(R^{\Lambda}(\beta)), \qquad
\xProj(R^{\Lambda}) := \bigoplus_{\beta \in Q^+}  \xProj(R^{\Lambda}(\beta)).
\]
\end{proposition}


\begin{corollary} \label{cor:QHSAblocks}
For any $\beta \in Q^+$,  $w\in W$, and reduced expression $w=r_{i_1} \dots r_{i_k}$ into simple transpositions, there is a superequivalence of supercategories
\[
\Theta_w \maps K^b(\xRep(R^{\Lambda}_{\beta})) \equiv
K^b(\xRep(R^{\Lambda}_{\Lambda-w(\Lambda-\beta)}))
\]
\end{corollary}

\begin{proof}
If $i\in I$ is odd, then the 2-supercategory $\mf{U}(\mf{g})$ contains a subcategory isomorphic to $\mf{U}(\mf{sl}_2)$ from Definition~\ref{def:oddU}, and if $i\in I$ is even, it contains a subcategory isomorphic to the usual categorification of $\mf{sl}_2$.  In both cases, there is a derived (super)equivalence $\Theta_i$ coming either from the usual Chuang-Rouquier complex in the even case, or the newly defined complexes from Definition~\ref{def:Theta}.  Hence, for any element $w$ in the Weyl group of type $\mf{g}$ and reduced expression $w=r_{i_1}\dots r_{i_k}$,  there is a derived superequivalence $\Theta_w := \Theta_{i_1} \dots \Theta_{i_k}$.
\end{proof}

\subsection{Abelian defect conjecture for spin symmetric groups}

In this section, let $\Bbbk$ be a field of odd characteristic $p=2\ell+1$.  Let $\mf{g}$ be the Kac-Moody Lie algebra of type $A_{2\ell}^{(2)}$ whose vertices of its Dynkin diagram are labelled $\{0,1,\dots, \ell\}$.  Let $\delta$ be the null root, $\{\alpha_i \mid i \in I\}$ the simple roots, $\{ \Lambda_i \mid i \in I\}$ the corresponding fundamental dominant weights, $Q_+$ the non-negative part of the root lattice, $P_+$ the set of dominant integral weights.

For $\beta =\sum_{i\in I}m_i\alpha_i\in Q_+$, let $R_{\beta}$ denote Kang-Kashiwara-Tsuchioka's~\cite{KKT} quiver-Hecke  superalgebra of type $\mf{g}$. For $\Lambda \in P_+$, we let $R_{\beta}^{\Lambda}$ denote the cyclotomic quiver Hecke superalgebra.  Likewise, let $RC_{\beta}$ and $RC_{\beta}^{\Lambda}$ denote the quiver Hecke-Clifford superalgebra.  Kang, Kashiwara, Tsuchioka show that there are super Morita equivalences
\[
RC_{\beta} \simeq_{sMor} R_{\beta} \otimes \cal{C}_{m_0}, \qquad
RC_{\beta}^{\Lambda} \simeq_{sMor} R_{\beta}^{\Lambda} \otimes \mf{C}_{m_0}
\]
where $\mf{C}_{m_0}$ is the Clifford superalgebra of size $m_0$.
They relate these cyclotomic quiver Hecke-Clifford superalgebras with blocks of the affine Sergeev superalgebra.  Our results then give new nontrivial derived superequivalences between blocks of the affine Sergeev superalgebra related by an action of the Weyl group $W$ of type $A_{2\ell}^{(2)}$.

When $\Lambda=\Lambda_0$, The Kang, Kashiwara, Oh categorification theorem~\cite{KKO} implies that $R^{\Lambda_0}_{\beta}$ is nonzero, if and only if $\Lambda_0-\beta$ is a nonzero weight of $V(\Lambda_0)$.  Nonzero weights of $V(\Lambda_0)$ can be understood in terms of the set of $p$-strict partitions $\lambda$ utilizing its residue content ${\rm cont}(\lambda)\in Q_+$, see \cite[Section 2.1d]{KLLi}.  Expressed in terms of the $\bar{p}$-core $\rho$ of the $p$-strict partition $\lambda$, the nonzero weights of $V(\Lambda_0)$ are of the form  ${\rm cont}(\rho)+d\delta$ for some $\bar{p}$-core partition $\rho$ and $d\in \Z_{\geq 0}$.  See \cite[Section 2.3a]{KLLi} and the references therein for details on the combinatorics of $p$-strict partitions and their cores.

Taking $\Lambda = \Lambda_0$, we have that, up to tensoring with a Clifford superalgebra, the cyclotomic quotient $R^{\Lambda_0}_{{\rm cont}(\rho) + d\delta}$ is   Morita superequivalent to a spin block $\cal{B}^{\rho,d}$ of the symmetric group.
%
Kleshchev and Livesey reduced the abelian defect conjecture for spin blocks to a conjecture~\cite[Conjecture 2]{KLLi} which implies the Kessar-Schaps conjecture~\cite{KeSc} and completes the program to prove Brou\'{e}'s abelian defect conjecture for the spin symmetric and alternating groups.  Our derived superequivalences give a proof of this conjecture.

\begin{theorem} [\cite{KLLi} Conjecture 2] \label{thm:spindef}
Let $\rho$ and $\rho'$ be $\bar{p}$-cores and $d$, $d' \in \Z_{\geq 0}$.  If $d=d'$, then the algebras $R^{\Lambda_0}_{{\rm cont}(\rho) + d\delta}$ and $R^{\Lambda_0}_{{\rm cont}(\rho') + d'\delta}$ are derived equivalent and $R^{\Lambda_0}_{{\rm cont}(\rho) + d\delta}\otimes \mf{C}_1$ and $R^{\Lambda_0}_{{\rm cont}(\rho') + d'\delta}\otimes \mf{C}_1$ are derived equivalent.
\end{theorem}

\begin{proof}

Assume that $d=d'$.  By the KKO categorification theorem, $R^{\Lambda_0}_{\beta}$ is nonzero if and only if $\Lambda_0-\beta$ is a nonzero weight of $V(\Lambda_0)$.   Furthermore, by \cite[Section 12]{Kac}, the nonzero weights of $V(\Lambda_0)$ are all of the form $w\Lambda_0 - d\delta$ with $w\in W$ and $d\in \Z_{\geq 0}$.
By \cite[Lemma 3.1.39]{KLLi} there is a map from the set of $p$-strict partitions $\cal{P}_p$ to the set of weights $P$ given by
\begin{align}
  \kappa \maps \cal{P}_p &\to P   \qquad
  \lambda  \mapsto \Lambda_0 - {\rm cont}(\lambda),
\end{align}
whose image is the set $\{w\Lambda_0 - d\delta \mid w \in W, d \in \Z_{\geq 0} \}$.  Further, this map restricts to a bijection between the set of $\bar{p}$-cores $\cal{C}_p$  and affine Weyl group orbits $W\Lambda_0$ of the highest weight $\Lambda_0$.
In particular, for $\bar{p}$-cores $\rho,\rho' \in \cal{C}_p$, we have
\[
 \Lambda_0-{\rm cont}(\rho) = w \Lambda_0, \qquad  \Lambda_0-{\rm cont}(\rho') = w' \Lambda_0
\]
for some affine Weyl group elements $w,w' \in W$.   Let $\tau = w'w^{-1} \in W$ so that $\tau w=w'$.
Write $\tau = r_{i_1} \dots r_{i_k}$   a reduced expression and let
\[
\Theta_{\tau} := \Theta_{i_1} \dots \Theta_{i_k}
\]
be the corresponding sequences of derived superequivalences defined in Section~\ref{sec:def-superequiv}.  Then by Corollary~\ref{cor:QHSAblocks}, $\Theta_{\tau}$ defines a derived superequivalence
from $\xRep(R^{\Lambda_0}_{{\rm cont}(\rho) + d\delta})$ to 
$\xRep(R^{\Lambda_0}_{\beta})$ where $\Lambda_0-\beta$ is the weight
\begin{align}
 \tau (\Lambda_0 -{\rm cont}(\rho)-d\delta )
    &= \tau(\Lambda_0 -{\rm cont}(\rho)) -\tau(d\delta) \nn
    = \tau(w\Lambda_0) - d\tau(\delta) \\
    &= w'\Lambda_0 - d\delta  = \Lambda_0 -{\rm cont}(\rho')-d\delta \nn
\end{align}
where we used that $\delta$ is $W$-invariant.  The theorem follows by \cite[Lemma 2.7]{KKT}.
\end{proof}

\begin{corollary}
Brou\'{e}'s abelian defect conjecture holds for double covers of the symmetric group and alternating groups.
\end{corollary}

\begin{proof}
The proof of this Corollary is reduced to Theorem~\ref{thm:spindef} in \cite[Theorem 5.4.12]{KLLi}.
\end{proof}

 \bibliographystyle{plain}
\bibliography{bib_spin}

\end{document}